\documentclass[11pt,twoside,a4paper]{article}
\usepackage{amsmath,amssymb,amsthm}
\newtheorem{Thm}{Theorem}[section]
\newtheorem{Lem}[Thm]{Lemma}
\newtheorem{Pro}[Thm]{Proposition}
\theoremstyle{definition}
\newtheorem{Def}[Thm]{Definition}

\numberwithin{equation}{section}

\newcommand{\R}{\mathbb{R}}
\newcommand{\Z}{\mathbb{Z}}

\newcommand{\cB}{\mathcal{B}}
\newcommand{\cC}{\mathcal{C}}
\newcommand{\cS}{\mathcal{S}}
\newcommand{\cU}{\mathcal{U}}

\newcommand{\diam}{\operatorname{diam}}
\newcommand{\inr}{\operatorname{int}}
\newcommand{\eps}{\epsilon}
\newcommand{\del}{\delta}
\newcommand{\Del}{\Delta}
\newcommand{\es}{\emptyset}
\newcommand{\gam}{\gamma}
\newcommand{\Lip}{\operatorname{Lip}}
\newcommand{\Ndim}{\dim_{\rm N}}
\newcommand{\asdim}{\operatorname{asdim}}
\newcommand{\lam}{\lambda}
\renewcommand{\phi}{\varphi}
\newcommand{\set}[2]{\{#1:\,#2\}}
\newcommand{\sm}{\setminus}
\newcommand{\sub}{\subset}
\newcommand{\st}{\operatorname{st}}
\newcommand{\tsum}{\textstyle\sum}
\newcommand{\sig}{\sigma}
\newcommand{\Sig}{\Sigma}

\begin{document}

\title{Nagata dimension, quasisymmetric embeddings, and Lipschitz extensions}
\author{Urs Lang \& Thilo Schlichenmaier
}
\markboth{U.~Lang \& T.~Schlichenmaier}{Nagata dimension}
\date{October 3, 2004}

\maketitle



\section{Introduction}

The large-scale geometry of metric spaces has been the subject of intense 
investigation in recent years. In his essay on asymptotic invariants of 
infinite groups~\cite{G2}, Gromov gave a comprehensive account of the 
area and introduced a number of new concepts. 
Here we discuss a variation of his notion of asymptotic 
dimension, with a view towards applications in analysis on 
metric spaces, an equally active field of current research.
The invariant considered was introduced and named Nagata dimension in 
a note by Assouad~\cite{As1}; indeed it is closely related to a theorem
of Nagata characterizing the topological dimension of metrizable spaces
(cf.~\cite[Thm.~5]{Nag1} or~\cite[p.~138]{Nag2}).
In contrast to the asymptotic dimension, the Nagata dimension of a 
metric space is in general not preserved under quasi-isometries, but it is
still a bi-Lipschitz invariant and, as it turns out, even a quasisymmetry 
invariant. The class of metric spaces with finite Nagata dimension includes
all doubling spaces, metric trees, euclidean buildings, and 
homogeneous or pinched negatively curved Hadamard manifolds, among others. 
One of our main results interposes between
theorems of Assouad~\cite{As2} and Dranishnikov~\cite{Dra}:
Every metric space with Nagata dimension at most $n$ admits in particular
a quasisymmetric embedding into the product of $n+1$ metric trees, 
cf.~Theorem~\ref{Thm:embeddingi} below. 
Another result asserts that a complete metric space with 
Nagata dimension $\le n$ is an absolute Lipschitz retract if and only if 
it is Lipschitz $m$-connected for $m = 0,1,\dots,n$; 
see Theorem~\ref{Thm:lip3i}.

We now give the precise definitions and outline the contents of the paper
in more detail. 
Suppose that $X = (X,d)$ is a metric space and $\cB = (B_i)_{i \in I}$ is a 
family of subsets of $X$.
The family $\cB$ is called \emph{$D$-bounded} for some constant $D \ge 0$ 
if $\diam B_i := \sup\set{d(x,x')}{x,x'\in B_i} \le D$
for all $i \in I$. The \emph{multiplicity} of the family is 
the infimum of all integers $n \ge 0$ such that each point in $X$ belongs to 
no more than $n$ members of $\cB$. For $s > 0$, the 
\emph{$s$-multiplicity} of $\cB$ is the infimum of all $n$ such that 
every subset of $X$ with diameter $\le s$ meets at most $n$ members of the 
family. The \emph{asymptotic dimension} $\asdim X$ of $X$ is defined as 
the infimum of all integers $n$ such that for all $s > 0$, 
$X$ possesses a $D$-bounded covering with $s$-multiplicity at most $n+1$ 
for some $D = D(s) < \infty$. This imposes no condition on small scales as 
it is not required that $D(s) \to 0$ for $s \to 0$.
This invariant was introduced (and denoted $\asdim_+ X$) 
by Gromov in~\cite[1.E]{G2}. He remarked that for various classes of 
spaces the function $D(s)$ can be chosen to be linear.

\begin{Def} \label{Def:ndim}
For a metric space $X$ the \emph{Nagata dimension} 
(or \emph{Assouad--Nagata dimension}) 
$\Ndim X$ of $X$ is the infimum of all integers $n$ with the 
following property:
There exists a constant $c > 0$ such that for all $s > 0$, 
$X$ has a $cs$-bounded covering with $s$-multiplicity at most $n+1$. 
\end{Def}

Note that this notion takes into account \emph{all} scales 
of the metric space in an equal manner. Clearly $\Ndim X \ge \asdim X$. 
The number $\Ndim X$ is unaffected if the covering sets are required to 
be open, or closed, or if the `test set' with diameter $\le s$ is replaced
by an open or closed ball of radius $s$ (the minimal constant $c$ may change, 
however). In particular,~\ref{Def:ndim} agrees with the definition 
given in~\cite{As1}.

In section~2 we gather a number of basic properties of the Nagata dimension.
It is easily seen that $\Ndim X \le \Ndim Y$ whenever $f \colon X \to Y$ is 
a map between metric spaces satisfying, for instance, 
$a\,d(x,x')^p \le d(f(x),f(x')) \le b\,d(x,x')^p$ 
for all $x,x' \in X$ and for some constants $a,b,p > 0$ (Lemma~\ref{Lem:fxy}). 
For every metric space $X$, the topological dimension $\dim X$ never 
exceeds $\Ndim X$ (Proposition~\ref{Pro:topdim}); 
as Assouad~\cite{As1} observed, the argument implicitly occurs 
in~\cite[p.~149]{Nag2}.
Each subset $X$ of $\R^n$ containing interior points satisfies 
$\Ndim X = n$. Every doubling metric space has finite Nagata dimension 
(Lemma~\ref{Lem:doubling}).
In Proposition~\ref{Pro:dimn} we characterize spaces with Nagata 
dimension $\le n$ in various ways.
We obtain the product formula $\Ndim(X \times Y) \le \Ndim X + \Ndim Y$ 
(Theorem~\ref{Thm:xxx}); the inequality may be strict. 
For $X = Y \cup Z$, the relation $\Ndim X = \sup\{\Ndim Y,\Ndim Z\}$ 
holds (Proposition~\ref{Pro:union}). 
Hence, every compact $n$-dimensional riemannian manifold 
$X$ satisfies $\Ndim X = n$.
 
In section~3 we determine the Nagata dimension, or prove that it is finite,
for certain hyperbolic or nonpositively curved spaces. 
This is again inspired by~\cite{G2} where the 
arguments are sketched for the asymptotic dimension.
Every product of $n$ non-trivial metric trees and every
euclidean building of rank $n$ has Nagata 
dimension $n$ (Propositions~\ref{Pro:trees} and~\ref{Pro:buildings}). 
By a \emph{metric tree} we mean a geodesic metric space all of whose geodesic 
triangles are degenerate, i.e.~isometric to tripods; no local finiteness
or compactness assumption is made. A geodesic metric space $X$ that is 
$\del$-hyperbolic in the sense of Gromov has finite Nagata dimension if
it satisfies the respective condition up to scale $\del$, i.e.~if there 
exist $n$ and $c$ such that for all $s \in (0,\del]$, 
$X$ has a $cs$-bounded covering with $s$-multiplicity at most $n+1$ 
(Proposition~\ref{Pro:hyperbolic}). 
Finally, every homogeneous Hadamard manifold has finite Nagata dimension
(Proposition~\ref{Pro:homogeneous}). 

Section~4 contains the proofs of Theorems~\ref{Thm:qsinvi} 
and~\ref{Thm:embeddingi} below.
A map $f$ from a metric space $X$ into another metric space $Y$
is called \emph{quasisymmetric} if it is injective and there exists a 
homeomorphism $\eta \colon [0,\infty) \to [0,\infty)$ such that
$d(x,z) \le t\,d(x',z)$ implies $d(f(x),f(z)) \le \eta(t)\,d(f(x'),f(z))$
for all $x,x',z \in X$ and $t > 0$. Then $f^{-1} \colon f(X) \to X$
is also quasisymmetric, and $f$ is a homeomorphism onto its image 
(see~\cite{TV} or~\cite[Chapt.~10]{He} for basic properties of 
quasisymmetric maps). 

\begin{Thm} \label{Thm:qsinvi}
Let $X,Y$ be two metric spaces, and let $f \colon X \to Y$ be a quasisymmetric
homeomorphism. Then $\Ndim X = \Ndim Y$.
\end{Thm}

Recall that Assouad's theorem~\cite{As2} asserts that for every doubling 
metric space $(X,d)$ and every exponent $p \in (0,1)$, there is an 
$N$ such that the metric space $(X,d^p)$ admits a bi-Lipschitz 
embedding into $\R^N$.
Dranishnikov~\cite{Dra} showed that every geodesic metric space with
bounded geometry and asymptotic dimension at most $n$ admits a large-scale
uniform embedding into the product of $n+1$ locally finite metric trees.

\begin{Thm} \label{Thm:embeddingi}
Let $(X,d)$ be a metric space with $\Ndim X \le n < \infty$. Then for
all sufficiently small exponents $p \in (0,1)$, there exists a bi-Lipschitz 
embedding of $(X,d^p)$ into the product of $n + 1$ metric trees.
\end{Thm}

In particular, $(X,d)$ admits a quasisymmetric embedding into the product of 
$n + 1$ metric trees.
In the general case,~\ref{Thm:embeddingi} is optimal 
with respect to the number of trees. On the other hand, 
the $n$-dimensional hyperbolic space admits a quasi-isometric embedding 
into the product of $n$ copies of a simplicial metric tree, see~\cite{BuS}.

Finally, in section~5, we resume our investigation~\cite{LS}, \cite{L},
\cite{LPS}, \cite{LPl} of the extendability of Lipschitz maps.  
We say that a pair of metric spaces $(X,Y)$ 
possesses the \emph{Lipschitz extension property} if there is a constant 
$C$ such that for every subset $Z \sub X$ and for every Lipschitz map 
$f \colon Z \to Y$, there is a Lipschitz extension 
$\bar f \colon X \to Y$ of $f$ with constant $\Lip(\bar f) \le C \Lip(f)$. 
A comprehensive characterization of such pairs is still missing. However,
we obtain complete results if one of the two spaces has finite Nagata 
dimension. We call a metric space $Y$ \emph{Lipschitz $m$-connected} for some
integer $m \ge 0$ if there is a constant $c_m$ such that every Lipschitz 
map $f \colon S^m \to Y$ has a Lipschitz extension 
$\bar f \colon B^{m+1} \to Y$ with constant $\Lip(\bar f) \le c_m\Lip(f)$;
here $S^m$ and $B^{m+1}$ denote the unit sphere and closed ball in $\R^{m+1}$
equipped with the induced metric. This condition is easily verified
in the presence of an appropriate weak convexity property of the metric. 
In particular, every Banach space and every geodesic metric space 
with convex metric (cf.~\cite[II.1.3]{BrH}) is Lipschitz $m$-connected 
for all $m \ge 0$.

\begin{Thm} \label{Thm:lip1i}
Suppose that $X,Y$ are metric spaces, $\Ndim X \le n < \infty$, and 
$Y$ is complete. If $Y$ is Lipschitz $m$-connected for $m = 0,1,\dots,n-1$,
then the pair $(X,Y)$ has the Lipschitz extension property.
\end{Thm}

As a corollary we obtain the fact that for a complete metric space 
$Y$, the pair $(\R^n,Y)$ has the Lipschitz extension 
property if and only if $Y$ is Lipschitz $m$-connected for $m = 0,\dots,n-1$.
This can be shown more directly by means of a Whitney cube decomposition of
the complement of a closed set $Z \sub \R^n$, cf.~\cite[Thm.~(1.2)]{Al} 
and~\cite[Thm.~2]{JLS}. The proof of Theorem~\ref{Thm:lip1i} may be viewed 
as a generalization of that argument. An application of~\ref{Thm:lip1i}
is the density of Lipschitz maps in various Sobolev classes of maps from
doubling metric measure spaces into Banach spaces or convex metric spaces, 
see for instance~\cite[Sect.~6]{HKST}.

\begin{Thm} \label{Thm:lip3i}
Suppose that $Y$ is a complete metric space with $\Ndim Y \le n < \infty$,
and $Y$ is Lipschitz $m$-connected for $m = 0,1,\dots,n$. Then $Y$ is an 
absolute Lipschitz retract; equivalently, the pair $(X,Y)$ has the Lipschitz 
extension property for every metric space $X$.
\end{Thm}

This result is obtained as a corollary of Theorem~\ref{Thm:lip2}
which provides Lipschitz extensions for maps $f \colon Z \to Y$ 
defined on a set $Z \sub X$ with $\Ndim Z \le n$. 
When combined with the results of section~3, Theorem~\ref{Thm:lip3i}
unifies and generalizes the results obtained 
in~\cite[1.2]{LPS} and~\cite[4.6 and~6.5]{LPl}.

\smallskip
{\it Acknowledgements.}
We thank Sergei Buyalo, Alexander Dranishnikov, Bruce Kleiner and 
Viktor Schroeder for inspiring discussions. Thanks to~\cite{DraZ} 
we became aware of~\cite{As1}. 


\section{Basic properties}

In this section we record a number of basic properties of the Nagata 
dimension defined in~\ref{Def:ndim}. 
We start with a simple lemma on the behavior under mappings. 

\begin{Lem} \label{Lem:fxy}
Suppose that $X,Y$ are metric spaces, $f$ maps $X$ into $Y$, and there
are two homeomorphisms $\phi,\psi \colon [0, \infty) \to [0, \infty)$
such that 
\[
\phi(d(x,x')) \le d(f(x),f(x')) \le \psi(d(x,x')) 
\quad \text{for all $x,x' \in X$}
\]
and $\sup_{s > 0} \frac{1}{s} \phi^{-1}(\bar c \psi(s)) < \infty$
for every constant $\bar c > 0$. Then $\Ndim X \le \Ndim Y$.
\end{Lem}

For instance, the condition on $\phi$ and $\psi$ is satisfied if 
$\phi(s) = a s^p$ and $\psi(s) = b s^p$ for all $s \ge 0$
and for some constants $a,b,p > 0$. 
Taking $p = 1$, we see that bi-Lipschitz homeomorphisms preserve $\Ndim$.

\begin{proof}
Suppose that $\Ndim Y = n < \infty$. Let $s > 0$. By definition there is a 
constant $\bar c > 0$, independent of $s$, and a $\bar c\psi(s)$-bounded 
covering $(C_i)_{i \in I}$ of $Y$ with $\psi(s)$-multiplicity at most 
$n+1$. Put $B_i := f^{-1}(C_i)$; $(B_i)_{i \in I}$ is
a covering of $X$ with $s$-multiplicity $\le n+1$. Whenever
$x,x' \in B_i$, then
\[
d(x,x') \le \phi^{-1}(d(f(x),f(x')) \le \phi^{-1}(\bar c\psi(s)) \le c s
\]
for some constant $c > 0$ independent of $s$. Hence $\Ndim X \le n$.
\end{proof}

Next we show that for every metric space $X$, $\Ndim X$ is at least as large
as the topological dimension $\dim X$ of $X$. For two families $\cU,\cU'$ of
subsets of $X$ we write $\cU' < \cU$ if each member of $\cU'$ is contained 
in some member of $\cU$. Recall that $\dim X$ equals the infimum of all 
integers $n$ such that for every open covering $\cU$ of $X$ there is
an open covering $\cU' < \cU$ of $X$ with multiplicity 
at most $n+1$.
A theorem due to~Vop\v{e}nka, cf.~\cite[p.~126]{Nag2}, states that 
$\dim X \le n$ if and only if there is a sequence of open coverings 
$\cU^1 > \cU^2 > \dots$ of $X$ such that each $\cU^k$ has 
multiplicity at most $n+1$ and is $D_k$-bounded, where $D_k \to 0$
as $k \to \infty$. 

We denote by $U(x,r)$ and $B(x,r)$ the open and closed ball, respectively,
with center $x$ and radius $r$ in $X$. Similarly, $U(A,r)$ will denote the 
open $r$-neighborhood of a set $A \sub X$.

\begin{Pro} \label{Pro:topdim}
For every metric space $X$, $\dim X \le \Ndim X$.
\end{Pro}

\begin{proof}
Suppose that $\Ndim X = n < \infty$. Let $s > 0$. Choose a $cs$-bounded 
covering $(B_i)_{i \in I}$ of $X$ with $s$-multiplicity at most 
$n + 1$. Denote by $U_i$ the open $s/2$-neighborhood of $B_i$;
then $\diam U_i \le (c+1)s$. The open covering 
$(U_i)_{i \in I}$ has multiplicity at most $n+1$, and
\[
(U(x,s/2))_{x \in X} < (U_i)_{i \in I} 
< (U(x,(c + 1)s))_{x \in X}.
\]   
Repeating this construction for $s_k := (2(c+1))^{-k}$, $k = 1,2,\dots$,
we find a sequence of open coverings $\cU^k$ satisfying
the assumption of Vop\v{e}nka's theorem mentioned above. 
Compare~\cite[p.~149]{Nag2}.
\end{proof}

The set $X = \{0,1,\frac12,\frac13,\dots\}$, equipped with the metric 
induced from $\R$, is an example of a compact metric space with 
$\dim X = 0$ and $\Ndim X = 1$.
Note also that for every metric space $X$, $\Ndim X$ equals
the Nagata dimension of the metric completion $\bar X$ of $X$. 

The following lemma shows in particular that every doubling metric space
(cf.~\cite[10.13]{He}) has finite Nagata dimension. 

\begin{Lem} \label{Lem:doubling}
Suppose that $X$ is a metric space, $n \ge 0$, and $s > 0$.
If for every $x \in X$ the closed ball $B(x,3s)$ can be covered by
$n+1$ sets of diameter $\le s$, then $X$ admits a covering $\cB$ by closed 
balls of radius $s$ such that $\cB = \bigcup_{k = 0}^n \cB_k$ and each 
family $\cB_k$ has $s$-multiplicity at most $1$; in particular, 
$\cB$ has $s$-multiplicity $\le n+1$.
\end{Lem}

\begin{proof}
Choose a maximal set $Z \sub X$ with the property that 
$d(z,z') > s$ whenever $z,z' \in Z$ and $z \ne z'$.
The family $\cB := (B(z,s))_{z \in Z}$ covers $X$. 
For every $z \in Z$, the ball $B(z,3s)$ can be covered by $n+1$ 
sets of diameter $\le s$, each of which contains no more than one element 
of $Z$. Thus $Z \cap B(z,3s)$ has cardinality at most $n+1$. 
Therefore there exists a coloring $\chi \colon Z \to \{0,1,\dots,n\}$ 
such that $\chi(z) \ne \chi(z')$ whenever $z,z' \in Z$ and 
$0 < d(z,z') \le 3s$ (cf.~\cite[2.4]{As2}). 
Then for every $k \in \{0,1,\dots,n\}$, the family $\cB_k$ of all balls 
$B(z,s)$ with $\chi(z) = k$ has $s$-multiplicity at most $1$. 
\end{proof}

Next we discuss various characterizations of the inequality $\Ndim X \le n$,
see Proposition~\ref{Pro:dimn} below. This is a variation of the discussion
in~\cite[1.E, 1.E$_1$]{G2} for the asymptotic dimension. Our rendition 
is designed to facilitate the proof of the product formula 
$\Ndim(X \times Y) \le \Ndim X + \Ndim Y$, 
cf.~Theorem~\ref{Thm:xxx}.

We use the following notion of a piecewise 
euclidean polyhedral complex. By a \emph{convex euclidean polyhedral cell} 
$C$ we mean a compact metric space isometric to the convex hull of a 
positive finite number of points in some euclidean space. 
Faces, edges, vertices, the interior $\inr C$, and the dimension of 
such a cell $C$ are defined in an obvious way.
An $m$-dimensional convex euclidean polyhedral cell with exactly $m+1$ 
vertices is called a \emph{euclidean simplex}.
By a \emph{piecewise euclidean polyhedral complex} 
$Z = (Z,\cC)$ we mean a metric space $Z$ together with a 
covering $\cC \sub 2^Z$ of $Z$ such that each 
$C \in \cC$ is a convex euclidean polyhedral cell, every
face of $C$ is in $\cC$, and any two elements of $\cC$ either are 
disjoint or intersect in a common face. Note that we do not require 
the metric on $Z$ to be intrinsic; on the other hand, the cells $C$ are 
already euclidean when equipped with the induced metric rather than the 
induced intrinsic metric. 
For $z \in Z$ we denote by $\st(z,Z)$ the \emph{open star} of $z$ in $Z$,
i.e.~the union of all $\inr C$ such that $z \in C \in \cC$.
The dimension of the complex $Z$ is the supremum of the dimensions 
of its cells. If each cell is a euclidean simplex,
then $Z$ is called a \emph{piecewise euclidean simplicial complex}.
Given a piecewise euclidean polyhedral complex $Z = (Z,\cC)$, 
we may pass to its \emph{first barycentric subdivision} $Z_1 = (Z,\cC_1)$ to 
obtain a piecewise euclidean simplicial complex isometric to $Z$.
  
For the proofs of~\ref{Pro:dimn} and~\ref{Thm:xxx} we need the following
lemma.
 
\begin{Lem} \label{Lem:lebesgue}
Let $Z$ be a piecewise euclidean simplicial complex of dimension 
$\le n < \infty$. 
Suppose there are constants $\del,\gamma > 0$ such that the 
following two properties hold for each pair of simplices $S,T$ of $Z$:
\begin{enumerate}
\item[\rm (i)]
If $S \cap T = \es$, then $d(S,T) \ge \del$.
\item[\rm (ii)]
If $S \cap T \ne \es$, then 
$d(x, S \cap T) \le \gamma d(x,y)$ for all $x \in S$ and $y \in T$.
\end{enumerate}
Then there is a constant $\lam > 0$ depending only on $\gamma$ and $n$ 
such that every open ball of radius $\lam \del$ in $Z$ is contained in 
the open star $\st(z,Z)$ of some vertex $z$ of $Z$.
\end{Lem}

\begin{proof}
Let $p \in Z$, and let $\lam$ be a positive constant, to be determined
later. Denote by $\cS$ the set of all simplices $S$ of $Z$ with 
$U(p,\lam\del) \cap \inr S \ne \es$. We show that $\bigcap \cS \ne \es$;
then $U(p,\lam\del) \sub \st(z,Z)$ for every vertex $z$ in $\bigcap \cS$.
Since $Z$ has dimension $\le n$, it suffices to prove that 
$\bigcap_{i = 1}^{n+2} S_i \ne \es$ whenever $S_1,\dots,S_{n+2} \in \cS$.
Given such $S_1,\dots,S_{n+2}$, pick $x_i \in B(p,\lam\del) \cap
S_i$ for every $i$, and let $k \in \{1,\dots,n+2\}$ be the 
maximal index such that $\bigcap_{i=1}^k S_i \ne \es$. 
We claim that for $j = 1,\dots,k$, there 
is point $y_j \in \bigcap_{i=1}^j S_i =: T^j$ such that
\[
d(p,y_j) \le \gamma_j \lam\del,
\]
where $\gamma_1 = 1$ and $\gamma_j = 1 + \gamma + \gamma\gamma_{j-1}$
for all integers $j \ge 2$. 
For $j = 1$, put $y_1 := x_1$; then $d(p,y_1) \le \lam\del 
= \gamma_1 \lam\del$. 
For $j = 2,\dots,k$, if $y_{j-1} \in T^{j-1}$ satisfies 
$d(p,y_{j-1}) \le \gamma_{j-1} \lam\del$, then condition~(ii) provides 
a point $y_j \in S_j \cap T^{j-1} = T^j$ such that 
$d(x_j,y_j) \le \gamma d(x_j,y_{j-1})$. Hence,
$d(p,y_j) \le d(p,x_j) + d(x_j,y_j) \le d(p,x_j) + \gamma d(x_j,y_{j-1}) 
\le (1 + \gamma)d(p,x_j) + \gamma d(p,y_{j-1}) \le \gamma_j \lam\del$.
This proves the claim. 
Now if $k < n + 2$, then we have $S_{k+1} \cap T^k = \es$, therefore
$d(S_{k+1},T^k) \ge \del$ by condition~(i).  
On the other hand, $d(p,y_k) \le \gamma_k \lam\del \le 
\gamma_{n+1} \lam\del$, thus $d(S_{k+1},T^k) \le d(x_{k+1},y_k) \le 
d(p,x_{k+1}) + d(p,y_k) < \lam\del + \gamma_{n+1} \lam\del$.
Choosing the constant $\lam > 0$ so that
$\lam(1 + \gamma_{n+1}) = 1$ we obtain a contradiction.
Hence $k = n + 2$. 
\end{proof}

\begin{Pro} \label{Pro:dimn}
Let $X$ be a metric space, and let $n \ge 0$ be an integer. Then the following 
properties are equivalent:
\begin{enumerate}
\item[\rm (1)]
$\Ndim X \le n$, i.e.~there exists a constant $c_1 > 0$ such that for 
all $s > 0$, $X$ has a $c_1s$-bounded covering 
with $s$-multiplicity at most $n+1$.
\item[\rm (2)]
There exists a constant $c_2 > 0$ such that for all $s > 0$, there exists
a $1$-Lipschitz map $f \colon X \to Y$ into some piecewise euclidean 
simplicial complex $Y$ of dimension $\le n$ such that every finite 
subcomplex of $Y$ is isometric to a subcomplex of some regular euclidean
simplex of edge length $s$, and $\diam f^{-1}(\st(y,Y)) \le c_2s$ for every 
vertex $y$ of $Y$.
\item[\rm (3)]
There exists a constant $c_3 > 0$ such that for all $s > 0$, there exists
a $1$-Lipschitz map $g \colon X \to Z$ into some piecewise euclidean
polyhedral complex $Z$ of dimension $\le n$ such that
every open ball of radius $s$ in $Z$ is contained in the open star
$\st(z,Z_1)$ of some vertex $z$ of the first barycentric subdivision $Z_1$, 
and $\diam g^{-1}(\st(z,Z_1)) \le c_3s$ for every vertex $z$ of $Z_1$.
\item[\rm (4)]
There exists a constant $c_4 > 0$ such that for all $s > 0$, $X$ admits 
a $c_4s$-bounded covering of the form $\cB = \bigcup_{k=0}^n \cB_k$ where 
each family $\cB_k$ has $s$-multiplicity at most $1$.
\end{enumerate}
\end{Pro}

In view of characterization~(3) we have $\Ndim \R^n \le n$. 
It follows that $\Ndim A = n$ whenever $A \sub \R^n$ contains interior 
points, for $n = \dim A \le \Ndim A \le \Ndim \R^n \le n$.

\begin{proof}
Clearly~(4) implies~(1) for $c_1 = c_4$.

We show that~(1) implies~(2). 
Let $r > 0$, and let $(B_i)_{i \in I}$ be 
a $c_1r$-bounded covering of $X$ with $r$-multiplicity at most $n+1$.
For each $i \in I$, define $\sig_i \colon X \to \R$ 
by $\sig_i(x) := \sup\{0, \frac{r}{2} - d(x,B_i)\}$. For every 
$x \in X$ we have $\sig_i(x) > 0$ for at most $n + 1$ indices $i$, 
and the sum $\bar\sig := \sum_{i \in I} \sig_i$ satisfies 
$\bar\sig \ge \frac{r}{2}$.
Consider the Hilbert space $\ell^2(I)$.
There exists a constant $\lam > 0$, depending only on $n$, 
such that the map $f \colon X \to \ell^2(I)$ defined by 
$f(x) := (\lam r\sig_i(x)/\bar\sig(x))_{i \in I}$
is $1$-Lipschitz. The image $f(X)$ lies in the $n$-skeleton $\Sig^{(n)}$ of 
the simplex $\Sig = \set{(y_i)_{i \in I}}{y_i \ge 0,\,
\sum_i y_i = \lam r}$ with edge length $s = \sqrt{2}\lam r$.
Denote by $Y$ the minimal subcomplex of $\Sig^{(n)}$ containing $f(X)$.
Let $y$ be the vertex of $Y$ corresponding to the index $i$.
If $x \in f^{-1}(\st(y,Y))$, then $\sig_i(x) > 0$,
hence $\diam f^{-1}(\st(y,Y)) \le \diam B_i + 2\frac{r}{2} 
\le (c_1 + 1)r$. This yields~(2).

Next we show that~(2) implies~(3). 
Let $r > 0$. By~(2) there exists a $1$-Lipschitz 
map $g \colon X \to Z$ into some piecewise euclidean 
simplicial complex $Z$ of dimension $\le n$ such that 
every finite subcomplex of $Z$ is isometric to a subcomplex of some regular
euclidean simplex of edge length $r$, and $\diam g^{-1}(\st(z,Z)) \le c_2r$ 
for every vertex $z$ of $Z$. It follows that there exist constants
$\eps > 0$ and $\gamma > 0$ depending only on $n$ such that the 
conditions of Lemma~\ref{Lem:lebesgue} are satisfied with $\del = \eps r$ 
and $\gamma$ for the first barycentric subdivision $Z_1$ of $Z$.
Hence, every open ball of radius $s = \lam \eps r$ in $Z$ is contained in 
the open star $\st(z,Z_1)$ of some vertex $z$ of $Z_1$, where the constant 
$\lam > 0$ depends only on $n$. Moreover, every such star $\st(z,Z_1)$ 
is contained in $\st(z',Z)$ for some vertex $z'$ of $Z$, thus
$\diam g^{-1}(\st(z,Z_1)) \le  \diam g^{-1}(\st(z',Z)) \le c_2 r$.

It remains to show that~(3) implies~(4). Let $s > 0$. Choose $r > s$,
and let $g \colon X \to Z$ be given as in~(3), with $s$ replaced by $r$.
For every vertex $z$ of $Z_1$, denote by $B_z$ the set of all $x \in X$
such that $U(g(x),r) \sub \st(z,Z_1)$. The sets $B_z$ cover $X$, and their 
diameter is $\le c_3r$. Each vertex $z$ of $Z_1$ is the barycenter 
of a unique cell $C_z$ of $Z$. For $k = 0,\dots,n$, define $\cB_k$ as the 
family of all $B_z$ such that $C_z$ has dimension $k$. 
Whenever $B_z,B_{z'} \in \cB_k$ for some $k$ and $z \ne z'$,
then $\st(z,Z_1) \cap \st(z',Z_1) = \es$. It follows that 
$d(x,x') \ge d(g(x),g(x')) \ge r > s$ for all $x \in B_z$ and $x' \in B_{z'}$,
so $\cB_k$ has $s$-multiplicity at most $1$.
\end{proof}

\begin{Thm} \label{Thm:xxx}
Let $X$ and $X'$ be two (non-empty) metric spaces. Then
$\Ndim(X \times X') \le \Ndim X + \Ndim X'$.
\end{Thm}

In general, the inequality may be strict, as the following simple
example shows. Equip $\Z$ and $I = [0,1]$ with the metric induced from $\R$.
Then $\Ndim \Z = \Ndim I = 1$, and it is easily checked that 
$\Ndim(\Z \times I) = 1$ as well.

\begin{proof}
Suppose that $\Ndim X = n < \infty$ and $\Ndim X' = n' < \infty$.
Let $r > 0$. According to~\ref{Pro:dimn}(2) there exists
a $1$-Lipschitz map $f \colon X \to Y$ into some piecewise euclidean 
simplicial complex $Y$ of dimension $\le n$ such that every finite 
subcomplex of $Y$ is isometric to a subcomplex of some regular euclidean 
simplex of edge length $r$, and $\diam f^{-1}(\st(y,Y)) \le c_2 r$ for every 
vertex $y$ of $Y$. There is a similar map $f' \colon X' \to Y'$ into some 
complex $Y'$ of dimension $\le n'$, where we assume the corresponding
diameter bound to hold for the same constant $c_2$ independent of $r$.
The product $Z := Y \times Y'$ is a piecewise euclidean polyhedral complex 
of dimension $\le n + n'$. 
From the properties of $Y$ and $Y'$ it follows that 
there exist constants $\eps > 0$ and $\gamma > 0$ depending only on $n$ and 
$n'$ such that the conditions of Lemma~\ref{Lem:lebesgue} are satisfied 
with $\del = \eps r$ and $\gamma$ for the first barycentric subdivision 
$Z_1$ of $Z$. Hence, every open ball of radius $s = \lam \eps r$ in $Z$ is 
contained in the open star $\st(z,Z_1)$ of some vertex $z$ of $Z_1$, where 
the constant $\lam > 0$ depends only on $n$ and $n'$. 
The product map $g := f \times f': X \times X' \to Y \times Y' = Z$ is 
1-Lipschitz. Every vertex $z$ of $Z_1$ is of the form $(y,y')$ for
some vertices $y$ of $Y_1$ and $y'$ of $Y'_1$; then $\st(z,Z_1) \sub 
\st(y,Y_1) \times \st(y',Y'_1)$. Hence, $\diam g^{-1}(\st(z,Z_1))
\le \sqrt{2}c_2 r$. The result follows from 
characterization~\ref{Pro:dimn}(3).
\end{proof}

\begin{Pro} \label{Pro:union}
Suppose that $X$ is a metric space, and $X = Y \cup Z$.
Then $\Ndim X = \sup\{\Ndim Y,\Ndim Z\}$.
\end{Pro}

\begin{proof}
Clearly $\sup\{\Ndim Y,\Ndim Z\} \le \Ndim X$.
To prove the other inequality, suppose that  $\Ndim Y,\Ndim Z \le n < \infty$.
There is a constant $c > 0$ such that for all $s > 0$,
$Z$ has a $cs$-bounded covering $(D_k)_{k \in K}$ with $s$-multiplicity at 
most $n + 1$, and $Y$ has a $c(3 + 2c)s$-bounded covering $(C_j)_{j \in J}$
with $(3 + 2c)s$-multiplicity at most $n + 1$. 
Given such collections, define a covering $(B_i)_{i \in I}$ of $X$ as follows. 
Assume that $J \cap K = \es$. Denote by $L$ the set of all indices $k \in K$
for which there is no pair of points $y' \in \bigcup_{j \in J}C_j$, 
$z' \in D_k$ with $d(y',z') \le s$. Put $I := J \cup L$. Choose maps 
$j \colon K \sm L \to J$ and $y,z \colon K \sm L \to X$ such that 
$y(k) \in C_{j(k)}$, $z(k) \in D_k$, and $d(y(k),z(k)) \le s$. 
Now define $B_i := C_i \cup \bigcup_{k \in j^{-1}\{i\}} D_k$ for $i \in J$
and $B_i := D_i$ for $i \in L$. Each $B_i$ has diameter at most
$c(3 + 2c)s + 2(1 + c)s$. 
Let $E \sub X$ be a set of diameter $\le s$.
If $E$ is disjoint from $\bigcup_{i \in J} C_i$, then $E$ meets at most 
$n+1$ members of $(D_k)_{k \in K}$, each of which belongs to exactly
one $B_i$. If $E \cap C_i \ne \es$ for some $i \in J$, then 
$k \in K \sm L$ whenever $D_k \cap E \ne \es$. Then the set 
$E \cup \set{y(k)}{D_k \cap E \ne \es}$ has diameter $\le (3 + 2c)s$
and, hence, meets no more than $n+1$ members of $(C_i)_{i \in J}$.
This shows that $(B_i)_{i \in I}$ has $s$-multiplicity at most $n+1$.
\end{proof}

As a consequence of Proposition~\ref{Pro:union}, every compact 
$n$-dimensional riemannian manifold $X$ satisfies $\Ndim X = n$. 
The following result shows in particular that all uniform 
(in a suitable sense) $n$-dimensional riemannian manifolds have 
`Nagata dimension $n$ in the small'.

\begin{Pro} \label{Pro:inthesmall}
Let $X$ be a metric space. Suppose there exist constants $c,r > 0$ and 
integers $n,N \ge 0$ such that for all $x \in X$, 
\begin{enumerate}
\item[\rm (i)]
the ball $B(x,r)$ has a $cs$-bounded covering with $s$-multiplicity 
at most $n+1$ for every $s \in (0,r]$, and 
\item[\rm (ii)]
the ball $B(x,3r)$ can be covered by $N+1$ sets of diameter $\le r$.
\end{enumerate}
Then there exist $c',r' > 0$ such that for every $s \in (0,r']$, 
$X$ has a $c's$-bounded covering with $s$-multiplicity at most $n+1$. 
\end{Pro}

\begin{proof}
By~(ii) and Lemma~\ref{Lem:doubling}, $X$ admits a covering $\cB$ by closed 
balls of radius $r$ such that $\cB = \bigcup_{k = 0}^N \cB_k$ and each 
family $\cB_k$ has $r$-multiplicity at most $1$.
For $k = 0,\dots,N$, put $Y_k := \bigcup\cB_k$. Using~(i) we see that 
for every $s \in (0,r]$, $Y_k$ has a $cs$-bounded covering with 
$s$-multiplicity at most $n+1$. Since $X = \bigcup_{k=0}^N Y_k$, by applying
the procedure from the proof of Proposition~\ref{Pro:union} $N$ times 
we obtain the result.
\end{proof}


\section{Buildings, hyperbolic and homogeneous spaces}

In this section we consider various classes of hyperbolic 
or nonpositively curved spaces. 
The following technical lemma is used in the proofs of~\ref{Pro:trees}
and~\ref{Pro:buildings} which determine the Nagata dimension for products
of metric trees and euclidean buildings.
 
\begin{Lem} \label{Lem:xy}
Let $X$ and $Y$ be metric spaces, and let $\lam,\mu > 0$. 
Suppose that there exist a $1$-Lipschitz map $f \colon X \to Y$ and a map 
$h \colon X \times [0,\infty) \to X$ with the following three properties:
\begin{enumerate}
\item[\rm (i)] 
Whenever $C \sub Y$ is a non-empty bounded set, there exists a point 
$y \in Y$ such that for all $x \in f^{-1}(C)$ 
there is an $x' \in f^{-1}\{y\}$ with $d(x,x') \le \lam \diam C$.   
\item[\rm (ii)] 
For all $x \in X$ and $t \ge 0$, $d(h(x,t),x) \le t$.
\item[\rm (iii)] 
If $x,x' \in X$, $f(x) = f(x')$, and 
$t \ge \mu d(x,x')$, then $h(x,t) = h(x',t)$.
\end{enumerate}
Then $\Ndim X \le \Ndim Y$.  
\end{Lem}

\begin{proof}
Suppose that $\Ndim Y = n < \infty$. There is a constant $c > 0$ such that
for all $s > 0$, $Y$ has a $cs$-bounded covering 
$(C_j)_{j \in J}$ with $s$-multiplicity at most $n+1$. 
Let such a collection of non-empty sets $C_j$ be given.
For every $j \in J$, put $X^j := f^{-1}(C_j)$ and pick a point 
$y_j \in Y$ and a map $\pi^j \colon X^j \to f^{-1}\{y_j\}$ such that 
\begin{equation}
d(x,\pi^j(x)) \le \lam cs \label{eq:pij}
\end{equation}
for all $x \in X^j$; see condition~(i). Since $f$ is 
$1$-Lipschitz, the covering $(X^j)_{j \in J}$ of $X$ has $s$-multiplicity
at most $n+1$. Hence, in order to prove $\Ndim X \le n$, it suffices 
to show that each $X^j$ can be partitioned into a $\tilde c s$-bounded 
collection $\cB^j$ with $s$-multiplicity at most $1$, for some $\tilde c > 0$.
Fix an index $j \in J$. Let $t := \mu(2\lam c + 1)s$, 
and define an equivalence relation on $X^j$ such that $x \sim x'$ if 
and only if $h(\pi^j(x),t) = h(\pi^j(x'),t)$. By~(\ref{eq:pij}) and 
condition~(ii), every equivalence class has diameter
at most $2(\lam cs + t) = 2(\lam c + \mu(2\lam c + 1))s =: \tilde c s$. 
Moreover, every set $D \sub X^j$ of diameter $\le s$ meets at most one 
such class. For if $x,x' \in D$, then $d(\pi^j(x),\pi^j(x')) \le 
2\lam cs + s$, thus $x \sim x'$ by condition~(iii) and the choice of $t$.   
\end{proof}

Recall that by a \emph{metric tree} $T$ we mean a geodesic metric space
such that every geodesic triangle $\Del$ in $T$ is degenerate, 
in the sense that each side of $\Del$ is contained in the union of the 
remaining two. Equivalently, $T$ is a $0$-hyperbolic geodesic metric space, 
cf.~(\ref{eq:hyperbolic}).

\begin{Pro} \label{Pro:trees}
Let $T_1,\dots,T_n$ be metric trees, each containing more than one point. 
Then $\Ndim(T_1 \times \ldots \times T_n) = n$.
\end{Pro}

\begin{proof}
Let $T$ be a metric tree.
To see that $\Ndim T \le 1$, choose a basepoint $z \in T$ and denote by
$f \colon T \to [0,\infty)$ the distance function to $z$.
Define $h \colon T \times [0,\infty) \to T$ such that
$h(x,t)$ is the point at distance $\min\{t,d(x,z)\}$ from $x$ on the 
geodesic segment $[x,z]$. Then the conditions
of Lemma~\ref{Lem:xy} are satisfied with $\lam = 1$ and $\mu = \frac{1}{2}$;
as for condition~(i), choose $y = \inf C$ for $C \sub [0,\infty)$.  
Hence $\Ndim T \le 1$.
Given non-trivial metric trees $T_1,\dots,T_n$, we conclude that
$\Ndim(T_1 \times \ldots \times T_n) \le n$ by Theorem~\ref{Thm:xxx} (product).
On the other hand, $T_1 \times \ldots \times T_n$ contains a set
isometric to a product of $n$ non-trivial intervals,
so $\Ndim(T_1 \times \ldots \times T_n) \ge n$.
\end{proof}

The proof of the following result relies partly on~\cite[Sect.~6]{LPS}; 
we refer to that paper and to~\cite{BrH} for information on 
euclidean buildings and nonpositively curved metric spaces. 

\begin{Pro} \label{Pro:buildings}
Let $X$ be a (simplicial or affine) euclidean building of rank $n$. Then
$\Ndim X = n$.
\end{Pro}

\begin{proof}
The building $X$ can be written as the union $\bigcup_{\nu \in N} A_\nu$
of a family of $n$-flats $A_\nu$, each asymptotic to a fixed chamber
$\Del \sub \partial_\infty X$ at infinity. There is a $1$-Lipschitz map 
$f \colon X \to \R^n$ such that $f|A_\nu \colon A_\nu \to \R^n$ is an 
isometry for every $\nu \in N$, cf.~\cite[6.1]{LPS}. Then
condition~(i) of Lemma~\ref{Lem:xy} is satisfied with $\lam = 1$;
given $C \sub \R^n$, pick $y \in C$ arbitrarily.  
Denote by $\xi \in \partial_\infty X$ the center of $\Del$, and define
$h \colon X \times [0,\infty) \to X$ such that $t \mapsto h(x,t)$ is the 
(unit speed) ray from $x$ to $\xi$. Then condition~(ii) 
of~\ref{Lem:xy} holds as well. To verify condition~(iii),
suppose that $x \in A_\nu$, $x' \in A_{\nu'}$, and $f(x) = f(x')$.  
Denote by $[x,\Del) \sub A_\nu$ the Weyl chamber with basepoint $x$ 
asymptotic to $\Del$, and put $a(t) := d(h(x',t),[x,\Del))$. 
There exist a constant $\mu > 0$, depending only on $X$, and a
$t_0 \ge 0$ such that $a(t) = 0$ for $t \ge t_0$ and
$a(t) \ge \frac{1}{\mu}(t_0 - t)$ for $t \in [0,t_0]$, 
cf.~\cite[6.2 and~(18)]{LPS}. In particular, $h(x',t) \in A_\nu$ 
for $t \ge \mu d(x,x')$ since $\mu d(x,x') \ge \mu a(0) \ge t_0$. 
As $f|A_\nu$ is an isometry, $t \mapsto f(h(x,t))$ is a ray 
in $\R^n$; similarly $t \mapsto f(h(x',t))$ is a ray. The two
are asymptotic and start at the same point $f(x) = f(x')$, 
so they coincide. Since $h(x',t) \in A_\nu$ for 
$t \ge \mu d(x,x')$, it follows that $h(x,t) = h(x',t)$ for 
$t \ge \mu d(x,x')$. Now Lemma~\ref{Lem:xy} shows that $\Ndim X \le
\Ndim \R^n = n$. Since each apartment $A_\nu \sub X$ is isometric
to $\R^n$, it is clear that $\Ndim X \ge n$.
\end{proof}

The next result is a variation of Lemma~\ref{Lem:xy}; it is used
in the proofs of~\ref{Pro:hyperbolic} and~\ref{Pro:homogeneous}
which deal with Gromov hyperbolic spaces and homogeneous Hadamard
manifolds. 

\begin{Lem} \label{Lem:xydel}
Let $X$ and $Y$ be metric spaces, and let $\lam,\mu,\del > 0$. 
Suppose that there exist a $1$-Lipschitz map $f: X \to Y$ and a map 
$h \colon X \times [0,\infty) \to X$ with the following three properties:
\begin{enumerate}
\item[\rm (i)] 
Whenever $C \sub Y$ is a non-empty bounded set,
there exists a point $y \in Y$ such that for all $x \in f^{-1}(C)$
there is an $x' \in f^{-1}\{y\}$ with $d(x,x') \le \lam \diam C$.   
\item[\rm (ii)] 
For all $x \in X$ and $t \ge 0$, $d(h(x,t),x) \le t$.
\item[\rm (iii)] 
If $x,x' \in X$, $f(x) = f(x')$, and 
$t \ge \mu d(x,x')$, then $d(h(x,t),h(x',t)) \le \del$.
\end{enumerate}
Suppose further that there exist $n \ge 0$ and $c > 0$ such that 
for all $s \in (0,\del]$, $X$ has a $cs$-bounded covering 
with $s$-multiplicity at most $n+1$.
If $\Ndim Y < \infty$, then $\Ndim X < \infty$.  
\end{Lem}

\begin{proof}
Suppose that $\Ndim Y = n_Y < \infty$.
In view of the assumption on $X$, it suffices to show that there 
exist $n_X \ge 0$ and $c_X > 0$ such that for all $s > \del$, 
$X$ has a $c_X s$-bounded covering with $s$-multiplicity at most $n_X + 1$.
Given $s > \del$, choose a $c_Y s$-bounded covering 
$(C_j)_{j \in J}$ of $Y$ with $s$-multiplicity at most $n_Y + 1$, 
$C_j \ne \es$. For every $j \in J$, put $X^j := f^{-1}(C_j)$ and pick a 
point $y_j \in Y$ and a map $\pi^j \colon X^j \to f^{-1}\{y_j\}$ such that 
\begin{equation}
d(x,\pi^j(x)) \le \lam c_Y s \label{eq:pijdel}
\end{equation}
for all $x \in X^j$; see condition~(i). 
Since $f$ is $1$-Lipschitz, the covering $(X^j)_{j \in J}$ of $X$ has 
$s$-multiplicity at most $n_Y + 1$. Hence, it suffices to show that for some
$c_X > 0$, each $X^j$ has a $c_X s$-bounded 
covering $\cB^j$ with $s$-multiplicity at most $n+1$;
then $\Ndim X \le n_X := (n + 1)(n_Y + 1) - 1$.
Fix an index $j \in J$. 
Choose a $c\del$-bounded covering $(D_k)_{k \in K}$ of $X$ 
with $\del$-multiplicity at most $n+1$.  
Let $t := \mu(2\lam c_Y + 1)s$, and define
$B^j_k := \set{x \in X^j}{h(\pi^j(x),t) \in D_k}$.
By~(\ref{eq:pijdel}) and condition~(ii),
each $B^j_k$ has diameter $\le 2\lam c_Y s + 2t + c\del \le 
(2\lam c_Y  + 2\mu(2\lam c_Y + 1) + c)s =: c_X s$. 
Let $E \sub X^j$ be a set of diameter $\le s$.
If $x,x' \in E$, then $d(\pi^j(x),\pi^j(x')) \le 
2\lam c_Y s + s$, thus $d(h(\pi^j(x),t),h(\pi^j(x'),t)) \le \del$ 
by condition~(iii) and the choice of $t$. This shows that 
the set $h(\pi^j(E),t)$ has diameter $\le \del$.
Since $(D_k)_{k \in K}$ has $\del$-multiplicity $\le n+1$, it 
follows that $\cB^j := (B^j_k)_{k \in K}$ has $s$-multiplicity at most $n+1$.
\end{proof}

We turn to Gromov hyperbolic spaces. For $\del \ge 0$, we call a 
geodesic metric space $X$ \emph{$\del$-hyperbolic} if 
\begin{equation}
d(y,y') \le \del \label{eq:hyperbolic}
\end{equation}
whenever $x,x',z \in X$, $y$ and $y'$ lie on geodesics from $z$ to $x$ and 
$x'$, respectively, and $d(y,z) = d(y',z) \le 
\frac12 (d(x,z) + d(x',z) - d(x,x'))$, cf.~\cite{G1}.

\begin{Pro} \label{Pro:hyperbolic}
Let $X$ be a $\del$-hyperbolic geodesic metric space for some $\del > 0$. 
Suppose that there exist $n \ge 0$ and $c > 0$ such that for all 
$s \in (0,\del]$, $X$ has a $cs$-bounded covering with $s$-multiplicity 
at most $n+1$. Then $\Ndim X < \infty$. 
\end{Pro}

Note that the condition involving $s$ is satisfied in particular
if $X$ is doubling up to some scale, cf.~Lemma~\ref{Lem:doubling}.
 
\begin{proof}
As in the proof of Proposition~\ref{Pro:trees},
choose a basepoint $z \in X$ and denote by
$f \colon X \to [0,\infty)$ the distance function to $z$.
Define $h \colon X \times [0,\infty) \to X$ such that
$h(x,t)$ is the point at distance $\min\{t,d(x,z)\}$ from $x$ on the 
geodesic segment $[x,z]$. Then the conditions
of Lemma~\ref{Lem:xydel} are satisfied with $\lam = 1$ and 
$\mu = \frac{1}{2}$. 
\end{proof}

Finally, we consider homogeneous Hadamard manifolds, i.e.,
complete simply connected riemannian manifolds with nonpositive 
sectional curvature and transitive isometry group.

\begin{Pro} \label{Pro:homogeneous}
Let $X$ be a homogeneous Hadamard manifold. Then $\Ndim X < \infty$.
\end{Pro}

\begin{proof}
The argument is very similar to the proof of 
Proposition~\ref{Pro:buildings}, where
Lemma~\ref{Lem:xydel} is used in place of~\ref{Lem:xy}.
Every homogeneous Hadamard manifold $X$ of algebraic rank $n \ge 1$ 
is foliated by a family $(A_\nu)_{\nu \in N}$ of $n$-flats $A_\nu$, 
each asymptotic to a fixed chamber $\Del \sub \partial_\infty X$ at infinity. 
Again there is a $1$-Lipschitz map $f \colon X \to \R^n$ such that 
$f|A_\nu \colon A_\nu \to \R^n$ is an isometry for every $\nu \in N$, 
cf.~\cite[6.1]{LPS}. Then condition~(i) of Lemma~\ref{Lem:xydel} is satisfied 
with $\lam = 1$.  
Denote by $\xi \in \partial_\infty X$ the center of $\Del$, and define
$h \colon X \times [0,\infty) \to X$ such that $t \mapsto h(x,t)$ is the 
ray from $x$ to $\xi$. 
Then condition~(ii) of~\ref{Lem:xydel} holds as well. 
Since $X$ is homogeneous, it follows from Proposition~\ref{Pro:inthesmall}
that there exist constants $\del,c > 0$ such that for all $s \in (0,\del]$,
$X$ has a $cs$-bounded covering with $s$-multiplicity at most $\dim X + 1$.
Hence the last assumption of~\ref{Lem:xydel} is satisfied.
Finally, to verify condition~(iii),
suppose that $x \in A_\nu$, $x' \in A_{\nu'}$, and $f(x) = f(x')$.  
Denote by $[x,\Del) \sub A_\nu$ the Weyl chamber with basepoint $x$ 
asymptotic to $\Del$, and put $a(t) := d(h(x',t),[x,\Del))$. 
There exist a constant $\mu > 0$, depending only on $X$ and $\del$, and 
a $t_0 \ge 0$ such that $a(t) \le \del/2$ for $t \ge t_0$ and
$a(t) \ge \frac{1}{\mu}(t_0 - t)$ for $t \in [0,t_0]$, 
cf.~\cite[6.3 and~(18)]{LPS}. In particular, $d(h(x',t),A_\nu) 
\le \del/2$ for $t \ge \mu d(x,x')$ since $\mu d(x,x') 
\ge \mu a(0) \ge t_0$. 
As $f|A_\nu$ is an isometry, $t \mapsto f(h(x,t))$ is a ray 
in $\R^n$; similarly $t \mapsto f(h(x',t))$ is a ray. The two
are asymptotic and start at the same point $f(x) = f(x')$, 
so they coincide. Let $t \ge \mu d(x,x')$. If $p_t$ denotes the point in
$A_\nu$ closest to $h(x',t)$, then $d(h(x,t),p_t) = d(f(h(x,t)),f(p_t)) 
= d(f(h(x',t)),f(p_t)) \le d(h(x',t),p_t) \le \del/2$, 
thus $d(h(x,t),h(x',t)) \le \del$. 
Hence $\Ndim X < \infty$ by Lemma~\ref{Lem:xydel}. 
\end{proof}


\section{Quasisymmetric embeddings}

In this section we prove Theorems~\ref{Thm:qsinvi} and~\ref{Thm:embeddingi}. 
We need the following iterated version of property~\ref{Pro:dimn}(4).
A corresponding result for the asymptotic dimension was shown 
in~\cite[Prop.~1]{Dra}.

\begin{Pro} \label{Pro:cr}
Suppose that $X$ is a metric space with $\Ndim X \le n < \infty$. Then there 
is a constant $c > 0$ such that for all sufficiently large $r > 1$,
there exists a sequence of coverings $\cB^j$ of $X$, $j \in \Z$, 
with the following four properties:
\begin{enumerate}
\item[\rm (i)] 
For every $j \in \Z$, we have $\cB^j = \bigcup_{k = 0}^n \cB^j_k$
where each $\cB^j_k$ is a $cr^j$-bounded family with $r^j$-multiplicity 
at most $1$.
\item[\rm (ii)]
For all $j \in \Z$ and $x \in X$, there exists a $C \in \cB^j$ that
contains the closed ball $B(x,r^j)$.
\item[\rm (iii)]
For every $k \in \{0,\dots,n\}$ and every bounded set $B \sub X$, there is
a $C \in \cB_k := \bigcup_{j \in \Z} \cB^j_k$ such that $B \sub C$.
\item[\rm (iv)]
Whenever $B \in \cB^i_k$ and $C \in \cB^j_k$ for some $k$ and $i < j$,
then either $B \sub C$ or $d(x,y) > r^i$ for every pair of points
$x \in B$, $y \in C$.
\end{enumerate}  
\end{Pro}

\begin{proof}
First we choose coverings $\cB^j$, $j \in \Z$, such that~(i)
holds with $r^j$ and $c$ replaced by $5r^j$ and $c'$, where $r > 1$ 
is to be specified below and $c'$ is the constant from 
Proposition~\ref{Pro:dimn}(4). We fix a basepoint $z \in X$ and assume 
without loss of generality that for all $m \in \Z$ and $k \in \{0,\dots,n\}$, 
there is a $C \in \cB^{m(n+1) + k}_k$ that contains $z$.
Replacing each $C \in \cB^j$, $j \in \Z$, by $\bigcup_{x \in C}B(x,r^j)$, 
without changing the notation, we obtain~(ii) and~(iii). 
After this modification, each family $\cB^j_k$ 
is $(5c' + 2)r^j$-bounded and has $3r^j$-multiplicity at most $1$.

We write $C \succ B$ if $B \in \cB^i_k$ and $C \in \cB^j_k$ for some 
$k$ and $i < j$ and there is a pair of points $x \in B$, $y \in C$
with $d(x,y) \le 3r^i$. 
Given $C \in \cB^j_k$, we denote by $\hat C$ the union of $C$ with 
all $B$'s for which there exists a chain
\[
C \succ C_1 \succ \ldots \succ C_{m-1} \succ C_m = B
\] 
for some $m \ge 1$ and $C_h \in \cB^{j_h}_k$, $h = 1,\dots,m$. 
In this situation it follows that $B \sub B(y,R(j_1))$ for some $y \in C$,
where
\[
R(j_1) := (5c' + 5)\tsum_{q = 0}^\infty r^{j_1-q} = 
\tfrac{5c' + 5}{r - 1} r^{j_1 + 1}.
\]  
We choose the constant $r$ initially so large that 
$R(j_1) \le r^{j_1 + 1}$. Note that $r^{j_1 + 1} \le r^j$.
Hence, for each $\cB^j_k$, the corresponding family
of all $\hat C$ with $C \in \cB^j_k$ is $(5c' + 4)r^j$-bounded and 
has $r^j$-multiplicity at most $1$.  
Now let $B \in \cB^i_k$ and $C \in \cB^j_k$ for some $k$ and $i < j$.
Suppose that there is a pair of points $x \in \hat B$, $y \in \hat C$ with
$d(x,y) \le r^i$. Since $B \in \cB^i_k$, we have that $d(x,x') \le r^i$ 
for some $x' \in B$. There is a chain 
$C = C_0 \succ C_1 \succ \ldots \succ C_m$ for some $m \ge 0$ and 
$C_h \in \cB^{j_h}_k$, $h = 0,\dots,m$, such that $y \in C_m$.
Let $l \in \{0,\dots,m\}$ be the largest index with
$j_l \ge i$; note that $j_0 = j > i$. Then there is a point $y' \in C_l$ 
such that $d(y,y') \le r^i$; for if $l < m$, then $j_{l + 1} + 1 \le i$.
Hence, $d(x',y') \le 3r^i$.
As $\cB^i_k$ has $3r^i$-multiplicity at most $1$, 
this implies that $i < j_l$ and $C_l \succ B$. In other words, 
$\hat B \sub \hat C$.
Replacing each $C$ by the corresponding $\hat C$ we obtain~(iv).
The result holds for $c = 5c' + 4$ and for all $r \ge 5c' + 6$.
\end{proof}

Now let $f$ be a quasisymmetric embedding from a metric space $X$ into 
another metric space $Y$. By definition, there is a homeomorphism 
$\eta \colon [0,\infty) \to [0,\infty)$ such that
\begin{equation}
d(x,z) \le t\,d(x',z) \quad \text{implies} \quad
d(f(x),f(z)) \le \eta(t)\,d(f(x'),f(z)) \label{eq:eta}
\end{equation}
for all $x,x',z \in X$ and $t \ge 0$. Denoting by $\bar\eta \colon 
[0,\infty) \to [0,\infty)$ the homeomorphism satisfying 
$\bar\eta(t) = 1/\eta^{-1}(1/t)$ for all $t > 0$, we obtain that
\begin{equation}
d(f(x),f(z)) \le t\,d(f(x'),f(z)) \quad \text{implies} \quad
d(x,z) \le \bar\eta(t)\,d(x',z) \label{eq:bareta}
\end{equation}
for all $x,x',z \in X$ and $t \ge 0$.
Quasisymmetric maps take bounded sets to bounded sets (cf.~\cite{TV} 
or~\cite[10.8]{He}); 
moreover, if $B,B' \sub X$ are two bounded sets 
with $B \cap B' \ne \es$, then 
\begin{equation}
\diam B \le t \diam B' \quad \text{implies} \quad
\diam f(B) \le 2\eta(2t) \diam f(B'). \label{eq:seta}
\end{equation}
To see this, suppose that $\diam B < t \diam B'$ for some $t > 0$.
Take $z \in B \cap B'$ and $x \in B$; then there is a point
$x' \in B'$ such that $d(x,z) \le 2t\,d(x',z)$. Hence
$d(f(x),f(z)) \le \eta(2t)\,d(f(x'),f(z)) \le \eta(2t) \diam f(B')$
by~(\ref{eq:eta}).
As this holds for arbitrary $x \in B$,~(\ref{eq:seta}) follows.
Analogously, by~(\ref{eq:bareta}), 
\begin{equation}
\diam f(B) \le t \diam f(B') \quad \text{implies} \quad
\diam B \le 2\bar\eta(2t) \diam B'. \label{eq:sbareta}
\end{equation}
Now we are in the position to prove the quasisymmetry invariance of the 
Nagata dimension.

\begin{proof}[Proof of Theorem~\ref{Thm:qsinvi}.]
Suppose that $\diam X > 0$ and $\Ndim X \le n < \infty$.
Let $\eta,\bar\eta \colon [0,\infty) \to [0,\infty)$ be homeomorphisms
such that~(\ref{eq:seta}) and~(\ref{eq:sbareta}) hold for the homeomorphism 
$f \colon X \to Y$.
Let numbers $c,r$ and coverings $\cB^j = \bigcup_{k = 0}^n \cB^j_k$ of 
$X$, $j \in \Z$, be given as in Proposition~\ref{Pro:cr}.
Choose constants $0 < \del \le 1$ and $\bar c > 0$ such that 
\[
2\bar\eta(2\del)c \le 1 \quad \text{and} \quad 2\eta(4cr) \le \del \bar c.
\]
Fix $s > 0$ such that $\bar c s < \diam Y \in (0,\infty]$; we show that 
$Y$ has a $\bar c s$-bounded covering with $s$-multiplicity at most $n+1$.
Using~\ref{Pro:cr}(ii) we see that for every $x \in X$ there is 
a maximal index $j(x) \in \Z$ such that there exist $k \in \{0,\dots,n\}$ and
$C \in \cB^{j(x)}_k$ with $B(x,r^{j(x)}) \sub C$ and 
$\diam f(C) \le \bar c s$; we pick such $k$ and $C$
and denote them by $k(x)$ and $C_x$. Then we choose a set $Z \sub X$ such 
that for every $x \in X$ there is a $z \in Z$ with $C_x \sub C_z$, and 
$C_z \not\sub C_{z'}$ whenever $z,z' \in Z$, $z \ne z'$.
(For the existence of such a set $Z$, note that every strictly
increasing sequence $C_{x_1} \sub C_{x_2} \sub \dots$ is finite;
this follows from~\ref{Pro:cr}(i) since the set $\{j(x_1),j(x_2),\dots\}$ 
is bounded.)
The family $(f(C_z))_{z \in Z}$ covers $Y$ and is $\bar c s$-bounded.
To prove that it has $s$-multiplicity at most $n+1$, we show that for every
$k \in \{0,\dots,n\}$, the subfamily consisting of all $f(C_z)$ 
with $k(z) = k$ has $s$-multiplicity at most $1$. By~\ref{Pro:cr}(i) and~(iv),
the members of this subfamily are pairwise disjoint.
Let $D \sub X$ be a set with $\diam f(D) \le s$; we assume that 
$Z_{k,D} := \set{z \in Z}{k(z) = k,\,C_z \cap D \ne \es}$ is non-empty and
that $D \not\sub C_z$ for all $z \in  Z_{k,D}$, in particular $\diam D > 0$. 
We must show that $Z_{k,D}$ consists of a single element.
Let $z \in Z_{k,D}$. By~\ref{Pro:cr}(ii) and the definition of $j(z)$, there 
is a $C' \in \cB^{j(z)+1}$ such that $B(z,r^{j(z)+1}) \sub C'$ and 
$\diam f(C') > \bar c s$. First we show that $\diam C_z \ge \frac12 r^{j(z)}$.
If this were not true, then the fact that 
$B(z,r^{j(z)}) \sub C_z \sub B(z,\frac12 r^{j(z)})$
and $D \not\sub C_z$ would imply $\diam D \ge \frac12 r^{j(z)}$,
thus $\diam C' \le cr^{j(z) + 1} \le 2cr \diam D$. Moreover, 
$D$ would meet $C' \supset B(z,\frac12 r^{j(z)}) \supset C_z$, 
therefore~(\ref{eq:seta}) would lead to the contradiction
$\bar c s < \diam f(C') \le 2\eta(4cr) \diam f(D) \le \del \bar c s 
\le \bar c s$. Hence, we have that 
$\diam C' \le cr^{j(z) + 1} \le 2cr\diam C_z$.
Since $z \in C_z \cap C'$, using~(\ref{eq:seta}) we infer that
\[
\bar c s \le \diam f(C') \le 2\eta(4cr) \diam f(C_z) 
\le \del \bar c \diam f(C_z),
\]
thus $\diam f(D) \le s \le \del \diam f(C_z)$. Now~(\ref{eq:sbareta}) yields
\[
\diam D \le 2\bar\eta(2\del) \diam C_z \le 2\bar\eta(2\del) cr^{j(z)} 
\le r^{j(z)}.
\]
Since $\diam D > 0$, this gives a lower bound on $j(z)$.
If we choose $z \in Z_{k,D}$ so that $j(z)$ is minimal, the
inequality $\diam D \le r^{j(z)}$ together with~\ref{Pro:cr}(i) and~(iv) 
implies that $Z_{k,D} = \{z\}$. 
\end{proof}

We proceed to the proof of the embedding theorem.

\begin{proof}[Proof of Theorem~\ref{Thm:embeddingi}.]
Let numbers $c,r$ and coverings $\cB^j = \bigcup_{k = 0}^n \cB^j_k$, 
$j \in \Z$, be given as in Proposition~\ref{Pro:cr}.
Now we write $B \prec C$ (or $C \succ B$) 
if $B \in \cB^i_k$ and $C \in \cB^j_k$
for some $k$ and $i < j$, $B \sub C$, and if there is no 
$B' \in \cB^{i'}_k$ with $i < i' < j$ and $B \sub B' \sub C$.
Note that if $B \in \cB^i_k$, $B_1 \in \cB^{i_1}_k$, $i \le i_1$,
$B \ne B_1$, and $B,B_1 \prec C$, then $d(B,B_1) \ge r^i$
due to~\ref{Pro:cr}(i) and~(iv). 

Let $p \in (0,1)$. Given $C \in \cB^j_k$, we denote by $d_{C,p}$ 
the largest pseudometric on $X$ satisfying $d_{C,p} \le d^p$ and 
$\sup\set{d_{C,p}(x,x')}{x,x' \in B} = 0$ for all $B \prec C$. 
More explicitly,
\begin{align}
d_{C,p}(x,x') = \inf\bigl\{ d(x,x')^p,\,\inf\bigl(d(x,B_1)^p 
  &+ \tsum_{h = 1}^{m-1}d(B_h,B_{h+1})^p \nonumber \\
  &+ d(x',B_m)^p \bigr) \bigr\}, \label{eq:dcp}
\end{align}
where the second infimum is taken over all $m \ge 1$ and all sequences
$B_1,B_2,\dots,B_m \prec C$. 
We claim that if $p$ is chosen sufficiently small, then
\begin{equation} \label{eq:dcpbound}
d_{C,p}(x,x') \ge (c+1)^{-p}\bigl(d(x,x') - cr^{j-1}\bigr)^p
\end{equation}
for all $x,x' \in X$ with $d(x,x') \ge cr^{j-1}$.
Suppose that $B \in \cB^i_k$, $B_1 \in \cB^{i_1}_k$, 
$B_2 \in \cB^{i_2}_k$, $i \le i_1,i_2$, $B_1 \ne B \ne B_2$,
and $B,B_1,B_2 \prec C$.  
Then $b_h := d(B,B_h) \ge r^i$ for $h = 1,2$.
The function $\phi(b_1,b_2) := b_1{}^p + b_2{}^p -
(b_1 + b_2 + cr^i)^p$ on $[r^i,\infty)^2$ achieves its infimum
at $(b_1,b_2) = (r^i,r^i)$; it is nonnegative if 
$p \le \log(2)/\log(2+c)$. Since $\diam B \le cr^i$, we infer that
\[
d(B_1,B_2)^p \le (b_1 + b_2 + cr^i)^p \le d(B_1,B)^p + d(B,B_2)^p
\]
for this choice of $p$.
This means that in~(\ref{eq:dcp}), it suffices to consider sequences 
$B_1,\dots,B_m \prec C$ such that the corresponding indices   
satisfy $i_1 \le \ldots \le i_l$ and $i_l \ge \ldots \ge i_m$ for some
$l \in \{1,\dots,m\}$. For every such sequence,
\[
d(x,x') \le d(x,B_l) + cr^{j-1} + d(x',B_l)
\]
since $\diam B_l \le cr^{i_l} \le cr^{j-1}$. Moreover,
\[
d(x,B_l) \le d(x,B_1) + \tsum_{h=1}^{l-1}(c+1)d(B_h,B_{h+1})
\]
since $\diam B_h \le cr^{i_h} \le cd(B_h,B_{h+1})$, and  
a similar estimate holds for $d(x',B_l)$. We conclude that
\begin{align*}
d(x,x') - cr^{j-1} \le (c+1)\bigl(d(x,B_1) 
&+ \tsum_{h=1}^{m-1}d(B_h,B_{h+1}) \\
&+ d(x',B_m)\bigr).
\end{align*}
Together with the subadditivity of the function $d \mapsto d^p$,
this gives~(\ref{eq:dcpbound}).

For every index $k \in \{0, \dots, n\}$ we define a metric tree $T_k$ 
as follows. For $C \in \cB^j_k$ we define $\tau_C \colon C \to \R$ by 
\begin{equation}
\tau_C(x) := \sup \bigl\{ 0, \inf\{d_{C,p}(x, X \sm C),
r^{p j}\} - r^{p(j-1)} \bigr\}, \label{eq:tauc}
\end{equation}
and we denote by $I_C$ an isometric copy of the possibly
degenerate segment $[0, \sup_{x \in C} \tau_C(x)]$.
The tree $T_k$ is then obtained by gluing $0 \in I_B$ to 
$\tau_C(B) \in I_C$ whenever $B,C \in \cB_k = \bigcup_{j \in \Z} \cB^j_k$ 
and $B \prec C$, where $\tau_C(B)$ denotes the unique value of $\tau_C$ 
on $B$. 
Indeed $T_k$ is connected due to~\ref{Pro:cr}(iii), and it contains no 
non-trivial loop by the definition of the relation $\prec$. 
For each $C \in \cB_k$, we have a canonical isometric inclusion 
$I_C \sub T_k$, and we denote by $z_C$ the point $0 \in I_C \sub T_k$.
Let $\bar T_k$ denote the metric completion of $T_k$. 
Given an infinite sequence $C_1 \succ C_2 \succ \ldots$ with 
$C_h \in \cB^{j_h}_k$ for all $h$, the sequence 
$z_{C_1},z_{C_2},\dots$ is Cauchy since
\begin{align*}
\tsum_{h = 1}^\infty d(z_{C_h},z_{C_{h+1}}) 
  &= \tsum_{h = 1}^\infty \tau_{C_h}(C_{h+1}) 
   \le \tsum_{h = 1}^\infty r^{p j_h} \\
  &\le \tsum_{q = 0}^\infty r^{p(j_1 - q)}
   = \tfrac{r^p}{r^p - 1} r^{p j_1}.
\end{align*}
We define a map $f_k \colon X \to \bar T_k$ as follows. 
Let $x \in X$. If there is an infinite sequence 
$C_1 \succ C_2 \succ \ldots$ with $x \in C_h \in \cB_k$ 
for all $h$, then $f_k(x)$ is defined as the limit point of the 
sequence $z_{C_1},z_{C_2},\dots$ in $\bar T_k$.
If there is no such sequence, then by~\ref{Pro:cr}(iii) there exists a 
minimal $j \in \Z$ such that $x \in C$ for some uniquely determined 
$C \in \cB_k^j$, and $f_k(x)$ is defined as the point 
$\tau_C(x) \in I_C \sub T_k$.

We show that $f_k$ is Lipschitz with respect to $d^p$ on $X$. 
Let $x,x' \in X$, $x \ne x'$. By~\ref{Pro:cr}(i) and~(iii) there is a 
minimal $j \in \Z$ such that $x,x' \in C$ for some $C \in \cB^j_k$.
We consider the case that there exist infinite sequences
$C \succ C_1 \succ C_2 \succ \ldots$ and 
$C \succ C'_1 \succ C'_2 \succ \ldots$ such that $x \in C_h \in \cB_k^{j_h}$ 
and $x' \in C'_h \in \cB_k^{j'_h}$ for all $h$; 
the other cases are similar. We assume without loss of generality that 
$j_1 \le j'_1$. Then $d(C_1,C'_1) \ge r^{j_1}$ by the choice of $j$ 
and by~\ref{Pro:cr}(i) and~(iv). We have that
\[
d(f_k(x), f_k(x')) 
  = d(f_k(x), z_{C_1}) + d(z_{C_1}, z_{C'_1}) + d(z_{C'_1}, f_k(x')).
\]
Using the triangle inequality for $d_{C,p}$ we see that
\begin{align*}
d(z_{C_1}, z_{C'_1}) 
&= |\tau_C(C_1) - \tau_C(C'_1)| \\  
&\le d_{C,p}(C_1, C'_1) 
\le d(C_1, C'_1)^p
\le d(x, x')^p.
\end{align*}
Moreover, since $r^{j_1} \le d(C_1,C'_1) \le d(x,x')$, 
\[
d(f_k(x), z_{C_1}) = \tsum_{h = 1}^\infty d(z_{C_h},z_{C_{h+1}}) 
  \le \tfrac{r^p}{r^p - 1} r^{p j_1}
  \le \tfrac{r^p}{r^p - 1} d(x,x')^p.
\]
Finally, if $d(f_k(x'), z_{C'_1}) > 0$, then there is a smallest
index $l \ge 1$ with $d(z_{C'_l},z_{C'_{l + 1}}) > 0$, and
\[
d(f_k(x'), z_{C'_1})
  = \tsum_{h = l}^\infty d(z_{C'_h},z_{C'_{h+1}})
  \le \tfrac{r^p}{r^p - 1} r^{p j'_l}.
\]
Since $\tau_{C'_l}(C'_{l + 1}) = d(z_{C'_l},z_{C'_{l + 1}}) > 0$, 
it follows from~(\ref{eq:tauc}) that
\[
r^{p(j'_l-1)} \le d_{C'_l,p}(C'_{l+1}, X \sm C'_l) \le 
d(C'_{l+1}, X \sm C'_l)^p \le d(x,x')^p.
\]
Hence,
\[
d(f_k(x'), z_{C'_1}) \le \tfrac{r^p}{r^p - 1} r^p \,d(x,x')^p.
\]
Combining these estimates we obtain that
\[
d(f_k(x), f_k(x')) 
\le \bigl(1 + \tfrac{r^p}{r^p - 1}(1 + r^p)\bigr) d(x,x')^p.
\]
In particular, the map 
\[
f := (f_0, f_1, \dots, f_n) \colon (X,d^p) \to 
\bar T_0 \times \bar T_1 \times \ldots \times \bar T_n
\]
is Lipschitz.

We show that $f$ is actually bi-Lipschitz. 
Let $x, x' \in X$. Choose $j \in\Z$ such that 
$cr^j < d(x, x') \le cr^{j + 1}$. 
By~\ref{Pro:cr}(ii) there exist a $k \in \{0,\dots,n\}$ and a $C \in \cB_k^j$ 
such that $U(x, r^j) \sub C$.
Since $\diam C \le cr^j < d(x, x')$ we have $x' \notin C$ and therefore 
\[
d(f(x),f(x')) \ge d(f_k(x), f_k(x')) \ge d(f_k(x),z_C) 
\ge \tau_C(x).
\]
Using~(\ref{eq:dcpbound}) and the fact that
$d(x,X \sm C) \ge r^j$ we infer that
\[
d_{C,p}(x, X \sm C) \ge (c+1)^{-p}\bigl(d(x, X \sm C) 
- cr^{j-1}\bigr)^p \ge \bigl(\tfrac{r-c}{c+1}\bigr)^p r^{p(j-1)}
\] 
for $r \ge c$. Recalling~(\ref{eq:tauc}) and the inequality
$d(x,x') \le cr^{j+1}$ we conclude that 
\begin{align*}
d(f(x),f(x')) 
  &\ge \tau_C(x) 
   \ge \bigl(\bigl(\tfrac{r-c}{c+1}\bigr)^p - 1\bigr) 
       r^{p(j-1)}\\
  &\ge c^{-p} r^{-2p} \bigl(\bigl(
       \tfrac{r-c}{c+1}\bigr)^p - 1\bigr) d(x,x')^p
\end{align*}
for $r$ sufficiently large.
\end{proof}


\section{Lipschitz extensions}

Finally, we prove the extension results for Lipschitz maps stated in the 
introduction.

\begin{proof}[Proof of Theorem~\ref{Thm:lip1i}.]
Let $Z \sub X$ be a closed set, and let $f \colon Z \to Y$ be a 
$\lam$-Lipschitz map.
Let $r \ge 2$ be a fixed number, to be specified below. For $i \in \Z$, define 
$R_i := \set{x \in X}{r^i \le d(x,Z) < r^{i+1}}$. 
Pick a set $N \sub X \sm Z$ that is maximal subject to the following condition:
Whenever $i \in \Z$ and $x,x' \in N \cap (R_i \cup R_{i+1})$, $x \ne x'$, 
then $d(x,x') \ge \frac14 r^i$. Let $N_i := N \cap R_i$. Note that the family 
of open neighborhoods $U(N_i,\frac14 r^i)$, $i \in \Z$, 
covers $X \sm Z$. Choose a retraction $\rho \colon Z \cup N \to Z$ such that 
$d(x,\rho(x)) \le r^{i+1}$ if $x \in N_i$.
For $z \in Z$, $x \in N_i$ and $x' \in N_j$, $j \ge i$, it follows that
\begin{equation}
d(\rho(x),z) \le d(x,z) + r^{i+1} \le (1 + r)\,d(x,z) \le 2r\,d(x,z)
\label{eq:rhoxz}
\end{equation}
since $d(x,z) \ge r^i$, and
\[
d(\rho(x),\rho(x')) 
\le d(x,x') + r^{i+1} + r^{j+1}
\le (1 + 4r + 4r^2)\,d(x,x')
\]
since $d(x,x') \ge \frac14 r^i$ if $i \le j \le i+1$ and
$d(x,x') \ge r^j - r^{j-1} \ge r^{j-1}$ if $j \ge i+2$.
Hence, $\rho$ is Lipschitz with constant $\Lip(\rho) \le 1 + 4r + 4r^2$.

Now we use the assumption $\Ndim X \le n$; let $c$ be the constant from
Definition~\ref{Def:ndim}.
For every $i \in \Z$ we pick a $2c r^i$-bounded covering 
$(D^i_l)_{l \in L_i}$ of $R_i$ with $2 r^i$-multiplicity at most $n+1$.
Let $C^i_l := U(D^i_l \cap N_i,\frac12 r^i)$; then
\begin{equation}
\diam C^i_l \le (2c + 1) r^i. \label{eq:diamcil}
\end{equation}
By choosing $r$ sufficiently large, depending on $c$, 
we arrange that
\begin{enumerate}
\item[(i)]
each family $(C^i_l)_{l \in L_i}$ has 
$2(2c + 1)r^{i-1}$-multiplicity at most $n+1$, and
\item[(ii)]
there is no triple of sets $C^{i-1}_j, C^i_k, C^{i+1}_l$
with $C^{i-1}_j \cap C^i_k \ne \es \ne 
C^i_k \cap C^{i+1}_l$.
\end{enumerate}
The respective conditions are 
$2(2c + 1)r^{i-1} + 2 \cdot \frac12 r^i \le 2 r^i$ and
$r^{i+1} - r^i \ge \frac12 r^{i-1} + (2c + 1) r^i + \frac12 r^{i + 1}$.
The family $\cC := (C^i_l)_{i \in \Z,\, l \in L_i}$ covers $X \sm Z$;
in fact the following stronger property holds. Define $\tau^i_l \colon
X \sm Z \to \R$ by 
\[
\tau^i_l(x) := \sup\bigl\{0, 2 - \tfrac{4}{r^i}d(x, D^i_l \cap N_i)\bigr\};
\] 
note that $C^i_l = \{\tau^i_l > 0\}$. 
Then $\sup\set{\tau^i_l}{i \in \Z,\, l \in L_i} \ge 1$ on $X \sm Z$.  

Clearly $\cC$ has multiplicity at most $2(n+1)$.
To obtain the sharp result, we must get rid of the factor $2$;
we proceed similarly as in the proof of Proposition~\ref{Pro:union}. 
Each family $(C^i_l)_{l \in L_i}$ is replaced by 
a new family $(B^i_k)_{k \in K_i}$ as follows. 
The new index set $K_i$ is the set of all $k \in L_i$ such that 
$C^i_k \ne \es$ and $C^i_k \cap C^{i+1}_l = \es$ for all $l \in L_{i+1}$.
For every $j \in L_{i-1} \sm K_{i-1}$ with $C^{i-1}_j \ne \es$, choose
an index $k_j \in L_i$ such that $C^{i-1}_j \cap C^i_{k_j} \ne \es$.
By~(ii), $k_j \in K_i$. For $k \in K_i$, define 
\[
B^i_k := C^i_k \cup \textstyle\bigcup_{k_j = k} C^{i-1}_j 
\quad \text{and} \quad
\sig^i_k := \sup\{\tau^i_k, \sup_{k_j = k} \tau^{i-1}_j\}.
\] 
Note that $B^i_k = \{\sig^i_k > 0\}$. Let $A$ be the set
of all pairs $(i,k)$ with $i \in \Z$ and $k \in K_i$.
The family $\cB := (B^i_k)_{(i,k) \in A}$ covers $X \sm Z$;
in fact $\sup\set{\sig^i_k}{(i,k) \in A} \ge 1$ on $X \sm Z$. 
We claim that $\cB$ has multiplicity at most $n+1$. Let $x \in X \sm Z$. 
Let $i$ be the maximal index such that $x \in C^i_l$ for some 
$l \in L_i$. For every $C^{i-1}_j$ containing $x$, there is a point 
$x_j \in C^{i-1}_j \cap C^i_{k_j}$ with $d(x,x_j) \le \diam C^{i-1}_j \le 
(2c + 1)r^{i-1}$. Now the claim follows from~(i).
 
Consider the Hilbert space $\ell^2(A)$. 
We know that $\bar\sig := \sum_{(i,k) \in A} \sig^i_k \ge 1$,
and we define $g \colon X \sm Z \to \ell^2(A)$ by 
\[
g(x) := (\sig^i_k(x)/\bar\sig(x))_{(i,k) \in A}.
\]
The image of $g$ lies in the $n$-skeleton $\Sig^{(n)}$ of the simplex
$\Sig := \set{(v_{(i,k)})_{(i,k) \in A}}{v_{(i,k)} \ge 0,\,
\sum_{(i,k) \in A} v_{(i,k)} = 1} \sub \ell^2(A)$. 
For every $(i,k) \in A$, we choose a point 
$x^i_k \in D^i_k \cap N_i \sub C^i_k \sub B^i_k$.
Let $h^{(0)} \colon \Sig^{(0)} \to Y$ be the map that sends 
the vertex $e_{(i,k)}$ of $\Sig$ to the point $f(\rho(x^i_k))$. 
For $m = 0,1,\dots,n-1$, we successively extend 
$h^{(m)}$ to a map $h^{(m+1)} \colon \Sig^{(m+1)} \to Y$
by means of the Lipschitz $m$-connectedness of $Y$. 
The resulting map $h := h^{(n)} \colon \Sig^{(n)} \to Y$ 
is Lipschitz on every closed simplex $S$ of $\Sig^{(n)}$,
with
\[
\Lip(h|S) \le C_1 \diam h(S^{(0)})
\]
for some constant $C_1$ depending on $n$ and the 
Lipschitz connected\-ness constants $c_0,c_1,\dots,c_{n-1}$.  
Finally, we define the extension $\bar f \colon X \to Y$ of $f$ such that 
\[
\bar f = h \circ g
\]
on $X \sm Z$.

It remains to show that $\bar f$ is Lipschitz.
Let $x \in R_q$, $q \in \Z$. 
Let $S$ be the minimal closed simplex of the complex
$\Sig^{(n)}$ containing $g(x)$. 
If $e_{(i,k)}$ is a vertex of $S$, then $\sig^i_k(x) > 0$ and
$i \le q + 1$. Using~(\ref{eq:diamcil}) and the fact that
$B^i_k = \{\sig^i_k > 0\}$ we see that 
\[
d(x,x^i_k) \le (2c + 1)(r^{i-1} + r^i) \le (4c+2)r^{q+1}.
\] 
Since $h(e_{(i,k)}) = f(\rho(x^i_k))$, it follows that 
\[
\Lip(h|S) \le C_1 \diam h(S^{(0)}) 
\le C_1(8c+4)\Lip(\rho)\lam r^{q+1} =: C_2 \lam r^{q+1}
\]
whenever $x \in R_q$ and $S$ is the minimal closed simplex of 
$\Sig^{(n)}$ containing $g(x)$. Thus, for $x^i_k$ as above,
\[ 
d(\bar f(x),f(\rho(x^i_k))) 
= d(h(g(x)),h(e_{(i,k)}))
\le 2\Lip(h|S) \le 2C_2 \lam r^{q+1}.
\]
Using~(\ref{eq:rhoxz}) we conclude that if $z \in Z$, then 
\begin{align*}
d(\bar f(x), \bar f(z)) 
&\le d(\bar f(x),f(\rho(x^i_k))) + d(f(\rho(x^i_k)),f(z)) \\
&\le 2C_2 \lam r^{q+1} + 2r\lam\, d(x^i_k,z) \\
&\le 2\bigl(C_2 + (4c+2)r\bigr)\lam r^{q+1} + 2r\lam\, d(x,z) \\
&\le 2\bigl(C_2 + (4c+2)r + 1\bigr)r\lam\, d(x,z) \\
&=: C_3\lam\, d(x,z) 
\end{align*}
since $d(x^i_k,z) \le d(x,z) + (4c+2)r^{q+1}$ and $d(x,z) \ge r^q$.

Finally, let $x \in R_q$ and $x' \in R_p$, $p \le q$. Suppose 
first that $d(x,x') \ge \frac14 r^{q-2}$. Choose points $z,z' \in Z$ with 
$d(x,z),d(x',z') \le r^{q+1}$. Then 
\begin{align*}
d(\bar f(x),\bar f(x')) 
&\le d(\bar f(x),\bar f(z)) + d(\bar f(x'),\bar f(z')) 
     + d(f(z),f(z')) \\
&\le C_3 \lam \bigl(d(x,z) + d(x',z')\bigr)
     + \lam\,d(z,z') \\
&\le (C_3 + 1)\lam\bigl(d(x,z) + d(x',z')\bigr) 
+ \lam\,d(x,x') \\
&\le \bigl((C_3 + 1)8r^3 + 1\bigr)\lam\, d(x,x')
\end{align*}
since $d(x,z) + d(x',z') \le 2r^{q+1} \le 8r^3\, d(x,x')$.
Now assume that $d(x,x') < \frac14 r^{q-2}$. Let $S,S'$ 
be the minimal closed simplices of $\Sig^{(n)}$ containing 
$g(x),g(x')$, respectively.
Choose $(i,k) \in A$ such that $\sig^i_k(x) \ge 1$;
then $\sig^i_k(x') > 0$ since $d(x,x') < \frac14 r^{q-2} \le \frac14 r^{i-1}$.
This means that $e_{(i,k)}$ is a common vertex of $S$ and $S'$,
in particular $S \cap S' \ne \es$. Then there is a point 
$v \in S \cap S'$ such that $d(g(x),v) + d(g(x'),v) \le 
\hat c \,d(g(x),g(x'))$ for some constant $\hat c$ 
depending only on $n$. It follows that
\begin{align*}
d(\bar f(x),\bar f(x')) 
&\le d(h(g(x)),h(v)) + d(h(g(x')),h(v)) \\
&\le C_2\lam r^{q+1} \bigl(d(g(x),v) + d(g(x'),v)\bigr) \\
&\le C_2\hat c\lam r^{q+1}\, d(g(x),g(x')).
\end{align*}
Using the fact that every $\sig^i_k$ with $\sig^i_k(x) > 0$ 
or $\sig^i_k(x') > 0$ is Lipschitz with constant 
$\Lip(\sig^i_k) \le 4r^{-(i-1)} \le 4r^{-(q-2)}$, we infer that 
$d(g(x),g(x')) \le \bar c r^{-(q-2)} d(x,x')$ for some
constant $\bar c$ depending only on $n$. Hence,
\[
d(\bar f(x),\bar f(x')) 
\le C_2\hat c\bar cr^3 \lam\, d(x,x').
\]
We conclude that $\bar f$ is a Lipschitz extension of $f$ with
$\Lip(\bar f) \le C \lam$ for some constant $C$ depending only
on $n$, $c$, and the Lipschitz connectedness constants $c_0,\dots,c_{n-1}$
of $Y$. 
\end{proof}

Theorem~\ref{Thm:lip3i} is in fact a special case of the following result. 

\begin{Thm} \label{Thm:lip2}
Suppose that $X,Y$ are metric spaces, $Z \sub X$ is a closed set with 
$\Ndim Z \le n < \infty$, and $f \colon Z \to Y$ is a Lipschitz map.
If $Y$ is Lipschitz $m$-connected for $m = 0,1,\dots,n$, then 
there is a Lipschitz extension $\bar f \colon X \to Y$ of $f$.
\end{Thm}

Unlike in Theorem~\ref{Thm:lip1i}, the assumption
on the Nagata dimension now refers to the domain of $f$ rather than 
$\bar f$, at the cost of the additional condition that $Y$ be Lipschitz
$n$-connected. The proof is similar to the foregoing and would be 
shorter if we assumed $Y$ to be Lipschitz $m$-connected for 
$m = 0,1,\dots,2n$.

\begin{proof}
Every metric space admits an isometric embedding into some Banach space.
Thus, for simplicity, we assume without loss of generality that $X$ is a 
geodesic metric space.
Suppose that $f \colon Z \to Y$ is $\lam$-Lipschitz.
For some fixed number $r \ge 3$, to be specified below, and for all $i \in \Z$,
define sets $R_i$, $N$, $N_i = N \cap R_i$ and a Lipschitz retraction
$\rho \colon Z \cup N \to Z$ exactly as in the first paragraph of the proof of 
Theorem~\ref{Thm:lip1i}. Recall that the family of open neighborhoods
$U(N_i,\frac14 r^i)$, $i \in \Z$, covers $X \sm Z$.

For every $i \in \Z$, pick a $3cr^{i+1}$-bounded covering 
$(D^i_l)_{l \in L_i}$ of $Z$ with $3r^{i+1}$-multiplicity at most $n+1$.
Let $C^i_l := U((\rho|N_i)^{-1}(D^i_l),\frac12 r^i)$ and
$K_i := \set{k \in L_i}{C^i_k \ne \es}$. Note that
\begin{equation}
\diam C^i_k \le (3cr + 2r + 1)r^i \le 3(c+1)r^{i+1}. \label{eq:diamcik}
\end{equation}
By choosing $r$ sufficiently large, depending on $c$, we arrange that 
each family $(C^i_k)_{k \in K_i}$ has 
$6(c+1)r^i$-multiplicity at most $n + 1$.
This is the case if $6(c+1)r^i + 
2 \cdot \tfrac12 r^i + 2 r^{i+1} \le 3r^{i+1}$.

We assume that $K_{i-1} \cap K_i = \es$. Let $J_{i-1}$ be the set of 
all $j \in K_{i-1}$ with $C^{i-1}_j \cap C^i_k = \es$ for all $k \in K_i$.
Choose a map 
\[
\alpha_i \colon K_{i-1} \to K_i \cup J_{i-1} =: A_i
\]
such that $C^{i-1}_j \cap C^i_{\alpha_i(j)} \ne \es$ for all 
$j \in K_{i-1} \sm J_{i-1}$ and $\alpha_i(j) = j$ for all $j \in J_{i-1}$.
Let
\[
B^i_k := C^i_k \cup \textstyle\bigcup_{\alpha_i(j) = k} C^{i-1}_j
\]
for $k \in K_i$ and $B^i_k := C^{i-1}_{\alpha_i(k)} = C^{i-1}_k$ 
for $k \in J_{i-1}$.
We claim that the family $(B^i_k)_{k \in A_i}$ has multiplicity at most $n+1$.
Let $x \in X \sm Z$. If $x \notin \bigcup_{l \in K_i} C^i_l$, then
$x$ belongs to at most $n+1$ members of $(C^{i-1}_j)_{j \in K_{i-1}}$,
each of which is contained in exactly one $B^i_k$. Now suppose that 
$x \in C^i_l$ for some $l \in K_i$. For every $C^{i-1}_j$ containing $x$, 
there is a point $x_j \in C^{i-1}_j \cap C^i_{\alpha_i(j)}$ with 
$d(x,x_j) \le \diam C^{i-1}_j \le 3(c+1)r^i$. 
Then it follows from the choice of $r$ that $x$ belongs to no more than
$n+1$ members of $(B^i_k)_{k \in A_i}$.

Let $t_i := \frac12 r^i$ and 
$T_i := \set{x \in X}{t_i \le d(x,Z) \le t_{i+1}}$ for $i \in \Z$.
Note that if $T_i \cap C^{i'}_k \ne \es$ for some $i' \in \Z$ and 
$k \in K_{i'}$, then $i' \in \{i-1,i\}$, because $r^{i-1} + \frac12 r^{i-2}
\le \frac32 r^{i-1} \le t_i$ and $r^{i+1} - \frac12 r^{i+1} \ge t_{i+1}$.
In particular, the family $(B^i_k)_{k \in A_i}$ covers $T_i$.
As in the proof of Theorem~\ref{Thm:lip1i}, we construct for every $B^i_k$ 
a Lipschitz function $\sig^i_k \colon X \sm Z \to [0,2]$ with 
$B^i_k = \{\sig^i_k > 0\}$ and $\Lip(\sig^i_k) \le 4r^{-(i-1)}$. Moreover, 
for all $x \in T_i$ there is a $k \in A_i$ such that $\sig^i_k(x) \ge 1$,
hence $\bar\sig^i(x) := \sum_{k \in A_i} \sig^i_k(x) \ge 1$.
For every $i \in \Z$, we obtain a Lipschitz map 
\[
g_i \colon T_i \to \ell^2(A_i), \quad 
g_i(x) := (\sig^i_k(x)/\bar\sig^i(x))_{k \in A_i},
\]
whose image lies in the $n$-skeleton $\Sig_i^{(n)}$ of the 
simplex $\Sig_i = \set{(v_k)_{k \in A_i}}{v_k \ge 0,\,
\sum_{k \in A_i} v_k = 1} \sub \ell^2(A_i)$.  
Then we define
\[
\bar g_i \colon T_i \to \Sig_i^{(n)} \times [0,1]
\]
such that $\bar g_i(x) = (g_i(x),\gam_i(x))$, where
$\gam_i(x) := (d(x,Z) - t_i)/(t_{i+1} - t_i)$. 
Note that $\gam_i$ is Lipschitz with constant $\Lip(\gam_i) 
\le 1/(t_{i+1} - t_i) \le r^{-i}$.
Let $\tilde\Sig_i$ be the subsimplex of $\Sig_i$ corresponding
to the subset $K_i \sub A_i$. The map 
$\alpha_{i+1} \colon K_i \to A_{i+1}$ defined above induces 
a canonical simplicial map
\[
\beta_{i+1} \colon 
\tilde\Sig_i^{(n)} \times \{1\} \to \Sig_{i+1}^{(n)} \times \{0\}. 
\]
We have that $\bar g_i(T_i \cap T_{i+1}) \sub \tilde\Sig_i^{(n)} \times \{1\}$.
For every $x \in T_i \cap T_{i+1}$, 
$\beta_{i+1}(\bar g_i(x)) = \bar g_{i+1}(x)$ 
since $\sig^i_k(x) = \sig^{i+1}_{\alpha_{i+1}(k)}(x)$ for all $k \in K_i$.

We construct a sequence of maps 
\[
h_i \colon \Sig_i^{(n)} \times [0,1] \to Y,
\]
$i \in \Z$, as follows.
In a first step, each $h_i$ is defined on $\Sig_i^{(n)} \times \{0\}$.
Choose points 
$x^i_k \in (\rho|N_i)^{-1}(D^i_k) \sub C^i_k \sub B^i_k$ 
for $k \in K_i$ and 
$x^i_k \in (\rho|N_{i-1})^{-1}(D^{i-1}_k) \sub C^{i-1}_k = B^i_k$
for $k \in J_{i-1}$. We define $h_i$ on $\Sig_i^{(0)} \times \{0\}$
so that $h_i(e_k,0) = f(\rho(x^i_k))$ for every vertex $e_k$ of $\Sig_i$.
Then we extend $h_i$ to $\Sig_i^{(n)} \times \{0\}$
by means of the Lipschitz $m$-connectedness of $Y$ for $m = 0,1,\dots,n-1$.
For the second step, we observe that for every $i \in \Z$ we already have 
a map $h_{i+1} \circ \beta_{i+1} \colon \tilde\Sig_i^{(n)} \times \{1\} \to Y$.
This gives $h_i$ on $\tilde\Sig_i^{(n)} \times \{1\}$.
In the last step, we put $h_i(e_k,1) := f(\rho(x^i_k))$
for every vertex $e_k \in \Sig_i^{(0)} \sm \tilde\Sig_i^{(0)}$, 
i.e.~$k \in A_i \sm K_i = J_{i-1}$. Then we extend $h_i$ to the whole
$(n+1)$-dimensional polyhedral complex $\Sig_i^{(n)} \times [0,1]$ by means 
of the Lipschitz $m$-connectedness of $Y$ for $m = 0,1,\dots,n$. 
For every closed simplex $S$ of the complex $\Sig_i^{(n)}$, 
\[
\Lip(h_i|S \times [0,1]) \le C_1 \diam h_i(S^{(0)} \times \{0,1\})
\]
for some constant $C_1$ depending on $n$ and the Lipschitz connectedness
constants $c_0,c_1,\dots,c_n$.
We define the extension $\bar f \colon X \to Y$ of $f$ such that
\[
\bar f|T_i = h_i \circ \bar g_i
\]
for all $i \in \Z$. This is consistent as 
$h_i \circ \bar g_i = h_{i+1} \circ \beta_{i+1} \circ \bar g_i 
= h_{i+1} \circ \bar g_{i+1}$ on $T_i \cap T_{i+1}$.  

It remains to show that $\bar f$ is Lipschitz. Since we assumed $X$ to be 
a geodesic space, it suffices to prove that all $\bar f|Z \cup T_i$ are 
Lipschitz with uniform constants.
Let $x \in T_i$, $i \in \Z$. 
Let $S$ be the minimal closed simplex of the complex
$\Sig_i^{(n)}$ containing $g_i(x)$. 
If $e_k$ is a vertex of $S$, then $\sig^i_k(x) > 0$.
Using~(\ref{eq:diamcik}) and the fact that $B^i_k = \{\sig^i_k > 0\}$
we see that 
\[
d(x,x^i_k) \le 3(c+1)(r^i + r^{i+1}) \le 4(c+1)r^{i+1}.
\]
Moreover, in case $k \in K_i$,
\[
d(x,x^{i+1}_{\alpha_{i+1}(k)}) \le 3(c+1)(r^i + r^{i+1} + r^{i+2}) 
\le 5(c+1)r^{i+2}.
\]
According to the definition of $h_i$ we have $h_i(e_k,0) = f(\rho(x^i_k))$, 
\[
h_i(e_k,1) = h_{i+1}(\beta_{i+1}(e_k,1))
= h_{i+1}(e_{\alpha_{i+1}(k)},0) = f(\rho(x^{i+1}_{\alpha_{i+1}(k)}))
\]
in case $k \in K_i$, and $h_i(e_k,1) = f(\rho(x^i_k))$
in case $k \in J_{i-1}$. It follows that 
\begin{align*}
\Lip(h_i|S \times [0,1]) &\le C_1 \diam h_i(S^{(0)} \times \{0,1\}) \\
&\le 10 C_1(c+1)\Lip(\rho)\lam r^{i+2} =: C_2 \lam r^{i+2}
\end{align*}
whenever $x \in T_i$ and $S$ is the minimal closed simplex of 
$\Sig^{(n)}$ containing $g_i(x)$. Thus, for $x^i_k$ as above,
\begin{align*} 
d(\bar f(x),f(\rho(x^i_k))) 
&= d(h_i(\bar g_i(x)),h_i(e_k,0)) \\
&\le 2\Lip(h_i|S \times [0,1]) \le 2C_2 \lam r^{i+2}.
\end{align*}
Using~(\ref{eq:rhoxz}) we conclude that if $z \in Z$, then 
\begin{align*}
d(\bar f(x),\bar f(z)) 
&\le d(\bar f(x),f(\rho(x^i_k))) + d(f(\rho(x^i_k)),f(z)) \\
&\le 2C_2 \lam r^{i+2} + 2r\lam\, d(x^i_k,z) \\
&\le 2\bigl(C_2 + 4(c+1)\bigr)\lam r^{i+2} + 2r\lam\, d(x,z) \\
&\le 2\bigl(2C_2r + 8(c+1)r + 1\bigr)r\lam\, d(x,z) \\
&=: C_3\lam\, d(x,z) 
\end{align*}
since $d(x^i_k,z) \le d(x,z) + 4(c+1)r^{i+1}$ and 
$d(x,z) \ge t_i = \frac12 r^i$.

Finally, let $x,x' \in T_i$. Suppose 
first that $d(x,x') \ge \frac14 r^{i-1}$. Choose points $z,z' \in Z$ with 
$d(x,z),d(x',z') \le t_{i+1} = \frac12 r^{i+1}$. Then 
\begin{align*}
d(\bar f(x),\bar f(x')) 
&\le d(\bar f(x),\bar f(z)) + d(\bar f(x'),\bar f(z')) 
     + d(f(z),f(z')) \\
&\le C_3 \lam \bigl(d(x,z) + d(x',z')\bigr)
     + \lam\,d(z,z') \\
&\le (C_3 + 1)\lam\bigl(d(x,z) + d(x',z')\bigr) 
+ \lam\,d(x,x') \\
&\le \bigl((C_3 + 1)4r^2 + 1\bigr)\lam\, d(x,x')
\end{align*}
since $d(x,z) + d(x',z') \le r^{i+1} \le 4r^2\, d(x,x')$.
Now assume that $d(x,x') < \frac14 r^{i-1}$. Let $S,S'$ 
be the minimal closed simplices of $\Sig_i^{(n)}$ containing 
$g_i(x),g_i(x')$, respectively.
Choose $k \in A_i$ such that $\sig^i_k(x) \ge 1$;
then $\sig^i_k(x') > 0$ since $d(x,x') < \frac14 r^{i-1}$.
This means that $e_k$ is a common vertex of $S$ and $S'$,
in particular $S \cap S' \ne \es$. Then there is a point 
$\bar v \in (S \cap S') \times [0,1]$ such that $d(\bar g_i(x),\bar v) 
+ d(\bar g_i(x'),\bar v) \le \hat c \,d(\bar g_i(x),\bar g_i(x'))$ 
for some constant $\hat c$ depending only on $n$. It follows that
\begin{align*}
d(\bar f(x),\bar f(x')) 
&\le d(h_i(\bar g_i(x)),h_i(\bar v)) + d(h_i(\bar g_i(x')),h_i(\bar v)) \\
&\le C_2\lam r^{i+2} \bigl(d(\bar g_i(x),\bar v) 
     + d(\bar g_i(x'),\bar v)\bigr) \\
&\le C_2\hat c\lam r^{i+2}\, d(\bar g_i(x),\bar g_i(x')).
\end{align*}
Since $\Lip(\sig^i_k) \le 4r^{-(i-1)}$ for all $k \in A_i$ 
and $\Lip(\gam_i) \le r^{-i}$, we see that 
$\bar g_i$ is $\bar c r^{-(i-1)}$-Lipschitz for some
constant $\bar c$ depending only on $n$. Hence,
\[
d(\bar f(x),\bar f(x')) 
\le C_2\hat c\bar cr^3 \lam\, d(x,x').
\]
We conclude that $\bar f$ is a Lipschitz extension of $f$ with
$\Lip(\bar f) \le C \lam$ for some constant $C$ depending only
on $n$, $c$, and the Lipschitz connectedness constants $c_0,\dots,c_n$
of $Y$. 
\end{proof}

\begin{proof}[Proof of Theorem~\ref{Thm:lip3i}.]
Whenever $Y$ sits isometrically in some metric space $X$, then the identity
map on $Y$ can be extended to a Lipschitz retraction from $X$ onto $Y$ 
according to Theorem~\ref{Thm:lip2}. 
Hence $Y$ is an absolute Lipschitz retract.
Equivalently, the pair $(X,Y)$ has the Lipschitz extension property for
every metric space $X$ (cf.~\cite[pp.~11ff]{BLi} for background on 
absolute Lipschitz retracts and various equivalent properties).
\end{proof}


\bigskip\noindent
Departement Mathematik\\
ETH Zentrum\\
R\"amistrasse 101\\
CH-8092 Z\"urich\\
Switzerland


\begin{thebibliography}{00}

\bibitem{Al} F. J. Almgren, The homotopy groups of the integral 
cycle groups, Topology 1 (1962), 257--299.

\bibitem{As1} P. Assouad, Sur la distance de Nagata,
C. R. Acad. Sci. Paris S\'er. I Math. 294 (1982), no. 1, 31--34.

\bibitem{As2} P. Assouad, Plongements lipschitziens dans $\R^n$, 
Bull. Soc. Math. France 111 (1983), 429--448.

\bibitem{BLi} Y. Benjamini, J. Lindenstrauss, Geometric Nonlinear 
Functional Analysis, Vol. 1, Amer. Math. Soc. Colloq. Publ., Vol. 48,
2000.

\bibitem{BrH} M. R. Bridson, A. Haefliger, Metric Spaces of 
Non-Positive Curvature, Springer 1999.

\bibitem{BuS} S. Buyalo, V. Schroeder, Embedding of hyperbolic spaces
in the product of trees, preprint 2003, math.GT/0311524.

\bibitem{Dra}
A. Dranishnikov, On hypersphericity of manifolds with finite asymptotic 
dimension, Trans. Amer. Math. Soc. 355 (2003), 155--167.

\bibitem{DraZ} A. Dranishnikov, M. Zarichnyi,
Universal spaces for asymptotic dimension, preprint 2002, math.GT/0211069.

\bibitem{G1} M. Gromov, Hyperbolic groups, pp.~75--263 in: 
S. M. Gersten (Ed.), Essays in Group Theory, Math. Sci. Res. Inst. Publ., 
Vol. 8, Springer 1987.

\bibitem{G2} M. Gromov, Asymptotic Invariants of Infinite Groups,
pp. 1--295 in: G. A. Niblo, M. A. Roller (eds.), Geometric Group Theory, 
Vol. 2, London Math. Soc. Lecture Note Series, no. 182, Cambridge Univ. 
Press 1993.

\bibitem{He} J. Heinonen, Lectures on Analysis on Metric Spaces,
Springer 2001.

\bibitem{HKST} 
J. Heinonen, P. Koskela, N. Shanmugalingam, J. T. Tyson,
Sobolev classes of Banach space-valued functions and quasiconformal 
mappings, J. Anal. Math. 85 (2001), 87--139.

\bibitem{JLS} W. B. Johnson, J. Lindenstrauss, G. Schechtman,
Extensions of Lipschitz maps into Banach spaces, Israel J. Math. 54 (1986),
129--138.

\bibitem{LS} U. Lang, V. Schroeder, Kirszbraun's theorem and metric spaces
of bounded curvature, Geom. Funct. Anal. (GAFA) 7 (1997), 535--560.

\bibitem{L} U. Lang, Extendability of large-scale Lipschitz maps,
Trans. Amer. Math. Soc. 351 (1999), 3975--3988.

\bibitem{LPS} U. Lang, B. Pavlovi\'{c}, V. Schroeder, 
Extensions of Lipschitz maps into Hadamard spaces, 
Geom. Funct. Anal. (GAFA) 10 (2000), 1527--1553.

\bibitem{LPl} U. Lang, C. Plaut, Bilipschitz embeddings of metric 
spaces into space forms, Geom. Dedicata 87 (2001), 285--307.

\bibitem{Nag1} J. Nagata, Note on dimension theory for metric spaces,  
Fund. Math. 45 (1958) 143--181.

\bibitem{Nag2} J. Nagata, Modern Dimension Theory, Noordhoff and 
North-Holland 1965.

\bibitem{TV} P. Tukia, J. V\"ais\"al\"a, 
Quasisymmetric embeddings of metric spaces,
Ann. Acad. Fenn. Ser. A I Math. 5 (1980), 97--114.

\end{thebibliography}
\end{document}